\pdfoutput=1
\documentclass{article}

\usepackage{makeidx}         
\usepackage{graphicx}        
\usepackage{multicol}        
\usepackage[bottom]{footmisc}
\usepackage{cite} 
 \usepackage{amssymb}
 \usepackage{amsmath}
 \usepackage{latexsym}

   \newtheorem{theorem}{Theorem}
   \newtheorem{proposition}{Proposition}
   \newtheorem{corollary}{Corollary}
  
  \newtheorem{lemma}{Lemma}
  \newtheorem{definition}{Definition}
  \numberwithin{equation}{section}

\newcommand{\qed}{\hspace*{\fill} $\Box$ \ifmmode \else
    \par\addvspace\topsep\fi}
\newenvironment {proof}{\par\addvspace\topsep\noindent {\it Proof.}
    \ignorespaces }

\makeindex             

\title{Graph Polynomials and Their Applications II: \\Interrelations
  and Interpretations} 
\author{Joanna A. Ellis-Monaghan\,$^1$
        \and
        Criel Merino\,$^2$
        }
\date{}

\begin{document}
 \maketitle

 \addtocounter{footnote}{1} \footnotetext{Department of Mathematics, Saint Michael's College, One
  Winooski Park, Colchester, VT, 05458, USA and Department of
  Mathematics and Statistics, University of 
  Vermont, 16 Colchester Avenue, Burlington, VT,  05405, USA.
\texttt{jellis-monaghan@smcvt.edu}}
 \addtocounter{footnote}{1} \footnotetext{Instituto de Matem\'aticas, Universidad Nacional Aut\'onoma
  de M\'exico, Area de la Investigaci\'on Cient\'{\i}fica,
Circuito Exterior, C.U., Coyoac\'an, 04510 M\'exico D.F., M\'exico. 
\texttt{merino@matem.unam.mx}}

\section{Introduction}
\label{sec:introduction}

 A graph polynomial is an algebraic object associated with a graph
 that is usually invariant at least under graph isomorphism.  As such, it
 encodes information about the graph, and enables algebraic methods
 for extracting this information.  This chapter surveys a 
 comprehensive, although not exhaustive, sampling of graph
 polynomials.  It concludes Graph Polynomials and their Applications I: The
 Tutte Polynomial by continuing the goal of
 providing a brief overview of a variety of techniques defining a
 graph polynomial and then for decoding the combinatorial information
 it contains.

The polynomials we discuss here are not generally specializations of the Tutte polynomial, but they are each in some way related to the Tutte polynomial, and often to one another.  We emphasize these
 interrelations and explore how an understanding of one polynomial can guide research into others.  We
 also discuss multivariable generalizations of some of these
 polynomials and the theory facilitated by this.  We conclude with two
 examples, one from biology and one from physics, that illustrate the
 applicability of graph polynomials in other fields.

\section{Formulating Graph Polynomials}
\label{sec:formulations}

We have seen two methods for formulating a graph polynomial with the
linear recursion (deletion/contraction) and generating function
definitions of the Tutte polynomial in the previous chapter.  Here we
will see several more.  We begin with one of the earliest graph
polynomials, the edge-difference polynomial, a multivariable
polynomial defined as a product and originally studied by
Sylvester~\cite{Syl78} and Peterson~\cite{Pet91} in the late 1800's.
More recently, it has been used to address list coloring questions
(see Alon and Tarsi~\cite{AT92} and Ellingham and Goddyn~\cite{EG96}),
where a list coloring of a graph is a proper coloring of the vertices
of a graph with the color of each vertex selected from a predetermined
list of colors assigned to that vertex.

\begin{definition}\label{edgediff}
\index{edge-difference polynomial}
\index{polynomial!edge-difference}
The edge-difference polynomial.
Let $\left( {v_1 , \ldots ,v_n } \right)$  be an ordering of the
vertices of a graph $G$.  Then $D\left( {G;x_1 , \ldots ,x_n } \right)
= \prod_{i < j} {( {x_i  - x_j } )} $, where the
product is over all edges $\left( {v_i ,v_j } \right)$  of G.
\end{definition}

Note that a proper coloring of $G$ corresponds to finding positive
integer values $N_i $  (not necessarily distinct) for each of the
$x_i$'s so that $D\left( {G;N_1 , \ldots ,N_n } \right)$ $\ne$ 0.

There are also several polynomials based on various determinants (or
even permanents; see Pathasarthy~\cite{Par89} for a survey) involving
the adjacency matrix of a graph.  Recall that $A(G)$, the adjacency
matrix of a graph, has entries $a_{ij}  = 1$ if $\left( {i,j} \right)$
is an edge of the graph and 0 if it is not. The characteristic
polynomial is the classic example of a such a graph polynomial, and
will be discussed further in Section~\ref{sec:polynomials}.

\begin{definition}\label{adjpolys} The characteristic polynomial. Let
  $A(G)$ be the adjacency matrix of a graph $G$. Then
  \index{characteristic polynomial}\index{polynomial!characteristic} 
         $f( {G;x} ) = |xI-A(G)|$.  
\end{definition}

 Other examples of such polynomials are the \emph{idiosyncratic polynomial}
 \index{idiosyncratic polynomial}\index{polynomial!idiosyncratic}
 introduced by   Tutte, see~\cite{Tut79}, that is defined by
 $\nu ( {G;x,y} ) = |A(G) + y(J-I-A) - xI |$, where $J$ is  the matrix
 having all  entries equal to 1. Also 
 $\mu (G;x,y ) = | xI-D(G)+A(G) |$ a polynomial introduced in
 \cite{Kel65},  where $D(G)$ is  the \emph{degree matrix} of  $G$ that
 is the diagonal matrix with  $\operatorname{deg}(i)$ in the position
  $(i,i)$. Note that $J-I-A$ is the adjacency matrix of the complement
 of $G$ and when   $G$ is a simple graph, $D(G)-A(G)$ is just
 the Laplacian matrix $L(G)$ of the  graph.

A number of important graph polynomials may be defined by state model
formulas. Loosely speaking, a state of a graph \index{graph state} is
some configuration resulting from making local assignments for
substructures (e.g. the edges or vertices) of the graph.  These
assignments may be, for example, associating an element of a given set
to each vertex, or even the result of reconfiguring the edges incident
with a vertex.  A graph polynomial is formed by associating an
expression, often a weighted monomial, to each state of the graph, and
then summing over all possible graph states. The language comes from
physics, and is also found in knot theory.  We will see several state
model graph polynomials among those surveyed below, as well as an
application of this method in the Potts model of statistical mechanics
in Section~\ref{sec:applications}.

\begin{figure}[hbtp]
      \begin{center}
       \includegraphics{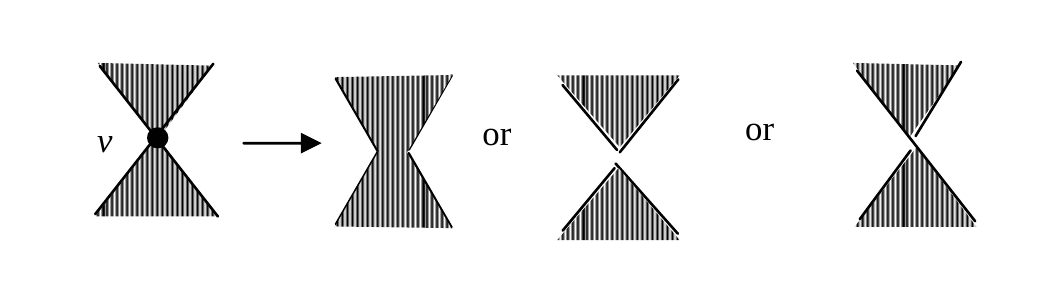}
        \caption{The three possible local reconfigurations at a vertex
        $v$, identified, from left to right, as white, black, and
        crossing.  Which strand passes over which in the crossing
        configuration does not affect the
        computation}\label{Fig:vertexfig}  
      \end{center}
      \end{figure}

An early example of a graph polynomial given by a state model
formulation \index{state model} is $P(G;x)$, the \emph{Penrose
polynomial}.\index{Penrose polynomial}\index{polynomial!Penrose}  This
polynomial graph invariant for planar graphs was defined implicitly by
Penrose~\cite{Pen69} in the context of tensor diagrams in physics, but
an excellent graph theoretical exposition can be found in
Aigner~\cite{Aig97}.  To compute $P(G;x)$, let $G$ be a plane graph,
and let $G_m $ be its medial graph, face two-colored with the
unbounded face colored white.  At each vertex, we consider three
possible local reconfigurations, as in Figure~\ref{Fig:vertexfig}.
A state $S$ of $G_m$ then results from choosing one of these three
reconfigurations at each vertex of $G_m$ and consists of a set of
disjoint closed curves (like a knot diagram).  Furthermore, to each
local reconfiguration at a vertex $v$, we assign a weight $\omega
(S,v)$ that is $+1$, 0, or $-1$ for a white, black, or crossing
configuration, respectively.

\begin{definition}\label{Penrose} The Penrose polynomial.  Let $G$ be
  a planar graph with medial graph $G_m$, and let $St(G_m)$ be the set
  of states of $G_m$ and let $St'(G_m)$ be the set of states with no
  black configurations.  Then,
\[
P(G;x) = \sum\limits_{St(G_m )}^{} {\left( {\left( {\prod\limits_{v
            \in G_m } {\omega (S,v)} } \right)x^{k(S)} } \right)}  =
         \sum\limits_{St'(G_m )}^{} 
            {\left( {\left( { - 1} \right)^{cr\left( S\right)}
            x^{k(S)} } \right)},  
\]
where $k(S)$ is the number of components in the graph state $S$, and
$cr(S)$ is
the number of crossing vertex configurations chosen in the state
$S$. 
\end{definition}

\begin{figure}[hbtp]
      \begin{center}
       \includegraphics[width=0.9\textwidth]{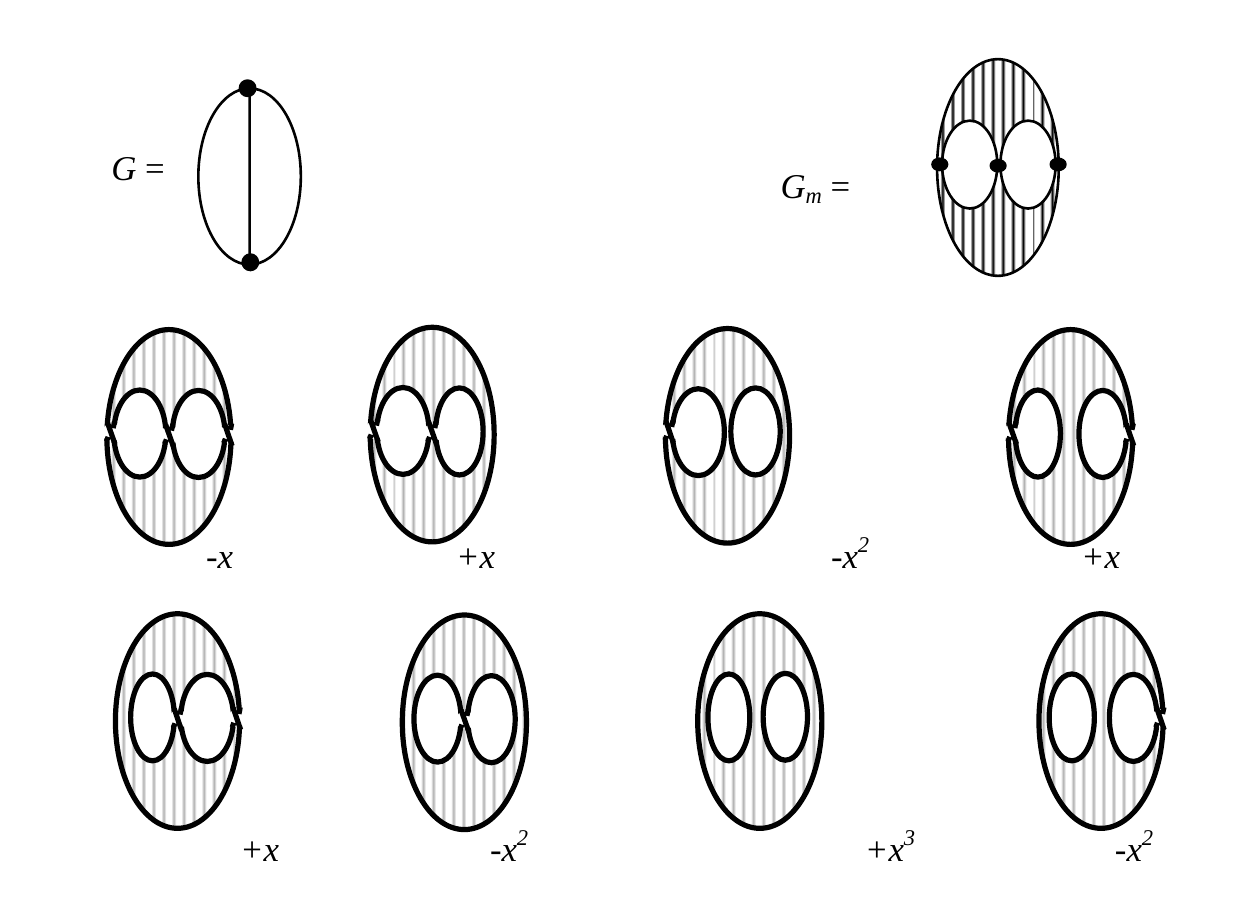}
        \caption{Computing the Penrose polynomial of a graph $G$ from
          the states of its medial graph}\label{Fig:Penrosefig}  
      \end{center}
      \end{figure}

For example, if $G$ is the $\theta$-graph consisting of two vertices
joined by three edges in parallel,
 then $P(G;x)=x^3-3x^2+2x$, as in~Figure~\ref{Fig:Penrosefig}. The
 Penrose polynomial may also be computed via a linear recursion
 relation (see Jaeger~\cite{Jae90} for example).  

The Penrose polynomial has some surprising properties, particularly
with respect to graph coloring.  The Four Color Theorem is equivalent
to showing that every planar, cubic, connected graph can be properly
edge-colored with three colors.  The Penrose polynomial, when applied
to planar, cubic, connected graphs, encodes exactly this information
(see Penrose~\cite{Pen69}):
\[
P\left( {G;3} \right) = \left( {\frac{{ - 1}}{4}}
\right)^{\frac{{\left| V \right|}}{2}} P\left( {G; - 2} \right) =
\text{the number of edge-3-colorings of } G.
\]

\section{Some Interrelated Polynomial Invariants}
\label{sec:polynomials} 
 
We present a further sampling of graph polynomials here.  They are each
related in some way to the Tutte polynomial, and have additional
relations among themselves.  These relations lead to combinatorial
insights as results for any one polynomial then inform those related
to it.  
 
 \subsection{Characteristic  and Matching Polynomials}
 \index{polynomial!characteristic}
 \index{polynomial!matching}
 
 The characteristic and matching polynomials are particularly
 interrelated, so we treat them together here, beginning with the
 characteristic polynomial $f(G;x)$ already introduced in
 Definition~\ref{adjpolys}.  Note that $f(G;x)$ is a monic polynomial
 of   degree $n$. Furthermore, since the adjacency matrix $A$ is real and
 symmetric, all its  eigenvalues are real, and thus all the zeros of
 $f(G;x)$ are real.  

 By using properties of
 determinants we can find interpretations of the coefficients of
 $f(G;x)$ in terms of the \emph{principal minors} of $A$. A principal
 minor of order $r$ is the determinant of an $r\times r$ submatrix of $A$
 obtained by choosing $r$ rows and  columns with the same set of
 indices. 
 \begin{proposition}\label{CharPolyMinors}
  Suppose that $f(G;x)=\sum_{i=0}^{n} a_{i}x^{n-i}$. Then $(-1)^{i}a_{i}$ is
  equal to the sum of the principal minors of $A$ with order $i$.
 \end{proposition}
 This property of the  characteristic polynomial can be found, for 
 example, in Horn and Johnson~\cite{HJ90}.
  
 Since the diagonal elements of $A$ are all zero, we have that
 $a_1=0$. The principal minors of order two and three which are not
 zero are of the form $|J-I|$, where $J$ is the matrix having all
 entries +1 and $J-I$ has order 2 or 3. The $2\times 2$ submatrices $J-I$
 of $A(G)$ correspond naturally to the edges of $G$ and the $3\times 3$
 submatrices $J-I$ correspond to the $K_3$ subgraphs of $G$. Thus
 $c_{2}=-|E(G)|$ and 
 $-c_{3}$ is  twice the number of $K_3$ subgraphs of $G$. 

 A \emph{linear
 subgraph}\index{linear subgraph}\label{linear_subgraph} of $G$ is a 
 subgraph whose components are edges or cycles.  An expression for the 
 coefficients of $f(G;x)$ in terms of linear subgraphs is given in the
 following.
  \begin{proposition}\label{CharPolyCoefficients}
 The coefficients of the characteristic polynomial may be expressed as
 \[
  (-1)^{i}a_{i} =
    \sum_{\Lambda} (-1)^{r(\Lambda)} 2^{r^{*}(\Lambda)} ,
\]
  where $r$ is the rank function and the sum is over all linear  subgraphs $\Lambda$ of $G$ having
  $i$ vertices. 
\end{proposition} 

 Note that because $\Lambda$ is a linear subgraph,
 $r^{*}(\Lambda)$ is simply the number of components in $\Lambda$
 that are cycles. The proof of Proposition~\ref{CharPolyCoefficients} 
 uses Proposition~\ref{CharPolyMinors} and 
 can be found in Harary~\cite{Har62}, while a detailed history of this result
is given by Cvetkovi\'c, Doob, and Sachs~\cite{CDS80}.

 As with the Tutte polynomial  we also have
 some reduction formulas and an expression for the derivative of the
 characteristic polynomial. 
\begin{theorem}\label{CharPolyProperties}
The characteristic polynomial of a graph satisfies the following
identities:
\begin{enumerate}
\item $f(G\cup H;x)=f(G;x)f(H;x)$, \label{CharPolyMult}
\item $f(G;x)= f(G\setminus e,x)- f(G-{u}-{v};x)$ if $e=\{u,v\}$ is
 a cut-edge  of $G$; \label{CharPolyLinear}
\item $\frac{\partial}{\partial x} f(G;x)=\sum_{v\in
    V(G)}f(G-{v};x)$. \label{CharPolyDerivative}
\end{enumerate}
\end{theorem} 

A proof of these properties can be found in
Godsil~\cite{God93}. Item~\ref{CharPolyMult} 
is an  easy exercises in matrix theory, as in Horn and
Johnson~\cite{HJ90}. Item~\ref{CharPolyLinear} can be proved by using
Proposition~\ref{CharPolyCoefficients} 
and considering the linear subgraphs of $G$ that use the edge $e$ and
the ones that do not use it. The result follows because $e$ is in no cycle of
$G$ if and only if it is a cut-edge. Item~\ref{CharPolyDerivative} can
be proved by using 
Proposition~\ref{CharPolyMinors} since any principal minor of order $i$ is
counted $n-i$ times in the right hand side of the formula in
Item~\ref{CharPolyDerivative}.  

 Given a graph $G$, the collection of  (unlabeled)
 subgraphs $G-v$, for $v\in V(G)$, is called the \emph{deck}
 \index{graph!deck} of
 $G$, and the individual subgraphs are called
 \emph{cards}.\index{graph!card} Thus, the deck for a graph on $n$ vertices
 consists of $n$ graphs, each 
 of which has $n-1$ vertices.   
 Ulam's reconstruction conjecture in \cite{Ula60} asserts that any finite graph
 $G$ with more than two vertices is uniquely determined by its deck,
 see~\cite{Ula60}, and we call any graph that satisfies the
 conjecture
 \emph{reconstructible}\index{graph!reconstructible}. Similarly, an
 invariant of $G$ which can be deduced from the 
 deck is called reconstructible.

 Clearly, the number $n$ of vertices of a graph is
 reconstructible. Also, as every edge is present in exactly $n-2$ cards,
 the number of edges is also reconstructible. The following useful
 example of reconstructibility is due to Kelly~\cite{Kel57}.   
\begin{lemma}[Kelly's lemma] Let $G$ and $H$ be graphs, and let
 $\nu(H,G)$ denote the number of 
 subgraphs of $G$ isomorphic to  $H$.  Then
 \[
   (|V(G)|-|V(H)|)\nu(H,G)=\sum_{v\in V(G)} \nu(H,G-v) .
\]
\end{lemma}
 The proof is by a double counting argument, and it follows that
 $\nu(H,G)$ is reconstructible whenever $|V(H)|<|V(G)|$. 

 Tutte  proved in~\cite{Tut67} that  the Tutte polynomial is
 reconstructible,  and thus  the chromatic polynomial, the flow
 polynomial,  the number of spanning trees and any invariant mentioned
 in the previous chapter are also reconstructible. Tutte  also proved that
 the characteristic polynomial is reconstructible in~\cite{Tut79}.  

\begin{theorem}\label{CharPolyReconstructible}
The characteristic polynomial of a graph is reconstructible.
\end{theorem} 

 For the proof, note that we have immediately from
 Theorem~\ref{CharPolyProperties} that $f'(G;x)$ is
 reconstructible. It then remains to prove that the constant term of 
 $f(G;x)$ is reconstructible. But by Proposition~\ref{CharPolyMinors},
 this is the same as proving that $|A(G)|$ is reconstructible. Then,
 using Theorem~\ref{CharPolyCoefficients} and an extension of Kelly's
 lemma (see~\cite{Koc81}), the problem is reduced to proving that the
 number of  Hamiltonian cycles is reconstructible.  A complete proof
 of Theorem~\ref{CharPolyReconstructible} based on the proof in
 Kocay~\cite{Koc81}, can be found  in~\cite{God93}. 

 Let us turn now to the matching polynomial.
 An $i$-\emph{matching}\index{matching@$i$-matching} in a graph $G$ is
 a set of $i$ 
 edges, no two of which 
 have a vertex in common. Let $\Phi_{i}(G)$ denote the number of
 $i$-matchings, and set $\Phi_{0}(G)=1$. Thus
 $\Phi_{1}(G)=m$ is the number of edges of $G$, and if $n$, the number
 of vertices, is even, then
 $\Phi_{n/2}(G)$ is the number of perfect matchings of $G$.  
 
\begin{definition} 
Let $G$ be a graph.  Then the
 matching polynomial\index{polynomial!matching}
 \index{matching polynomial}  of $G$ is
\[
   \mu(G;x)=\sum_{i\geq 0} (-1)^{i} \Phi_{i}(G) x^{n-2i} .
\] 
\end{definition} 

A more natural polynomial might be the \emph{matching generating
 polynomial},\index{matching generating polynomial}
 \index{polynomial!matching generating} given as the
 generating function of 
 $i$-matchings by
\[
    g(G;x)=\sum_{i\geq 0} \Phi_{i}(G) x^{i} .
\]
 However, the two polynomials are related by the identity
\[
      \mu(G;x)=x^{n}g(G;(-x)^{2}) ,
\]
  so there is no essential difference between them. 

 The matching polynomial is also known as the \emph{acyclic
 polynomial}\index{polynomial!acyclic} in Gutman and
 Trinajsti\'c~\cite{GT76}, 
 \emph{matching defect polynomial}\index{polynomial!matching defect}
 in Lov\'asz and Plummer~\cite{LP86} and \emph{reference 
 polynomial}\index{polynomial!reference} in Aihara~\cite{Aih76}. It
 has appeared independently in several different contexts. 
 In combinatorics, it was probably introduced by Farrell in~\cite{Far79a},
 but since the matching polynomial is essentially the same as the rook
 polynomial\index{polynomial!rook} for bipartite graphs (see
 Farrell~\cite{Far88}),  then its 
 origin can be traced back  at least to Riordan~\cite{Rio58}. In statistical
 physics it  appears because of  the monomer-dimer problem and was
 introduced by 
 Heilmann and Lieb in~\cite{HL70} and independently by Kunz
 in~\cite{Kun70}.  Finally, in theoretical chemistry was introduced by Hosaya
 in~\cite{Hos71} and later in connection with the so-called
 topological resonance energy by Gutman, Milun and Trinajsti\'c
 in~\cite{GT76,GMT76,GMT77} and independently by Aihara
 in~\cite{Aih76}. For a  full account of the history of the matching
 polynomial 
 see Gutman~\cite{Gut91}. 

 As with the Tutte and characteristic polynomials, we have some
 reduction formulas for the matching polynomial. The proof of the
 following theorem can be found in Godsil~\cite{God93}.
\begin{theorem}
 The matching polynomial satisfies the following identities:
 \begin{enumerate}
 \item $\mu(G\cup H;x)=\mu(G;x)\mu(H;x)$,\label{MatchPolyDisjoint}
 \item $\mu(G;x)=\mu(G\setminus e;x)-\mu(G-{u}-{v};x)$ if $e=\{u,v\}$ is
 an edge  of $G$,\label{MatchPolyDeletionEdge}
\item $\mu(G;x)=x\mu(G-{u};x)-\sum_{\{u,v\}\in E(G)}
  \mu(G-{v}-{u};x)$, if $u\in V(G)$, \label{MatchPolyDeletionVertex}
\item $\frac{\partial}{\partial x} \mu(G;x)=\sum_{v\in
    V(G)}\mu(G-{v};x)$. \label{MatchPolyDerivative} 
 \end{enumerate}
\end{theorem}
 To choose an $i$-matching in $G\cup H$ you need to choose an
 $s$-matching in $G$ and a $t$-matching in $H$ such that
 $s+t=i$. Item~\ref{MatchPolyDisjoint} then
  follows by the fundamental counting
 principle. For Item~\ref{MatchPolyDeletionEdge}, notice that the set of
 $i$-matchings can be partitioned into those $i$-matchings that use the
 edge $e$ and those that do not use
 it. Item~\ref{MatchPolyDeletionVertex} follows similarly. Finally,
 every $i$-matching of 
 $G$ with $i<n/2$ is counted $n-2i$ times in the right-hand side of
 the formula in Item~\ref{MatchPolyDerivative}, so the result follows.

 When the graph $G$ is a forest, a linear subgraph of $G$ with $j$ vertices  
 corresponds to a matching covering $j$ vertices, with
 $j$ even. Thus, Proposition~\ref{CharPolyCoefficients} has the
 following corollary, observed by Hosaya~\cite{Hos71}  and by Heilmann
 and Lieb~\cite{HL72}. 
 
\begin{corollary}\label{MatchPolyCharPoly}
 If $G$ is a forest then $f(G;x)=\mu(G;x)$.
\end{corollary}

  An  unexpected property of the matching polynomial, proved by
 Heilmann and Lieb~\cite{HL72},  is that all its zeros 
 are real, and furthermore the zeros for any graph $G$ interlace with
 the zeros of any of the cards in its deck.  The same paper also gives
 bounds for the zeros.  

 \begin{theorem} \label{MatchPolyzeros}
  For any graph $G$, the matching polynomial $\mu(G;x)$ has only real
  zeros. Furthermore, if 
  $u$ is any vertex in $G$ and if $a_1,a_2,\ldots, a_n$ are the zeros
  of  $\mu(G;x)$ while the zeros of $\mu(G-{u};x)$ are
  $a'_1,a'_2,\ldots, a'_{n-1}$, then
 \[ 
    a_1\leq a'_1\leq a_2\leq a'_2\leq \ldots\leq a'_{n-1}\leq a_{n} ,
 \]
 that is, the zeros of $\mu(G;x)$ and $\mu(G-{u};x)$ interlace. 
 \end{theorem}

 \begin{theorem}
 The (real) zeros, $a_1,a_2,\ldots, a_n$, of $\mu(G;x)$,  satisfy
\[        |a_i|<  2\sqrt{maxdeg(G)-1}. \]

 \end{theorem}

We outline a proof of  
 Theorem~\ref{MatchPolyzeros} from Godsil~\cite{God81b} that uses some of the
 results already 
 mentioned for $\mu(G;x)$ and $f(G;x)$. First, given a graph $G$ and a
 vertex $u$ in $G$, the 
 \emph{path tree}\index{path tree} $T(G,u)$ is the tree that has as
 its vertices the paths in $G$ 
 which start at $u$, and where two such vertices are
 joined by an edge if one represents a maximal proper subpath
 (i.e. all but the last edge) of the other. We then have 
 the following proposition from~\cite{God81b} that leads to a proof of
 Theorem~\ref{MatchPolyzeros}. 

\begin{proposition} Let $u$ be a vertex in a graph $G$, let $T=T(G,u)$ be
 the path tree of $G$ with respect to $u$, and let $u^{\prime}$ be the
 vertex of $T$ corresponding to the path of length 0 beginning at
 $u$. Then 
\[
   \frac{\mu(G-u;x)}{\mu(G;x)} = \frac{\mu(T-u^{\prime};x)}{\mu(T;x)} =
   \frac{f(T-u^{\prime};x)}{f(T;x)} .   
\]
\end{proposition}
 The last equality follows from
 Corollary~\ref{MatchPolyCharPoly}. Because all the roots of the
 characteristic  polynomial are real, we conclude that all zeros and
 poles of the  rational function ${\mu(G-u;x)}/{\mu(G;x)}$ 
 are real. An induction argument on the number
 of vertices in $G$ then yields the conclusion that all the zeros of
 $\mu(G;x)$ 
 are real.

 An  interesting combinatorial consequence of
 Theorem~\ref{MatchPolyzeros} is the following result of Heilmann and
 Lieb~\cite{HL72}, which gives a stark contrast with how little is known
 about the coefficients of the chromatic polynomial.

 \begin{theorem}
 For any graph $G$, the sequence $\Phi_{0}(G)$, $\Phi_{1}(G)$,
 $\ldots$ of coefficients of $g(G;x)$ is log-concave, that is
 $\Phi_i^{2}\geq \Phi_{i-1}\Phi_{i+1}$. 
 \end{theorem}

  The characteristic polynomial has been well studied, particularly with
  respect to graphs with the same characteristic polynomial.
  Godsil~\cite{God93} gives a thorough treatment of both the
  characteristic and 
  the matching polynomials. Another  good reference for the
  characteristic polynomial is Biggs~\cite{Big96}. Just as the matching
  polynomial is a way to study the matchings of a graph, the
  characteristic polynomial is a way to study the spectra of the
  adjacency matrix of a graph. Cvetkovi\'c, Doob, and Sachs have
  written a book~\cite{CDS80} dedicated to 
  the spectra of the adjacency matrix, and Lov\'asz and
  Plummer~\cite{LP86} have a book 
  devoted to the theory of matchings.   Furthermore, although  the
  characteristic 
  polynomial is not a complete invariant of graphs, it is conjectured
  that the characteristic polynomial of a graph $G$ is reconstructible
  from its polynomial deck, i.e. from the set 
  of characteristic polynomials of the cards of $G$.  See Gutman and
  Cvetkovi\'c~\cite{GC75} for the conjecture, 
  and then Cvetkovi\'c and Lepovi\'c~\cite{CL98}, where it is proved in
  the case of trees.

 \subsection{Ehrhart Polynomial}

 A \emph{convex polytope}\index{polytope!convex} $P$ is the convex
 hull of a finite set of  
 points in ${\mathbb  R}^m$. We denote the interior
 of $P$ (in the usual topological sense) by $P^o$. A convex polytope $P$ is said to be a 
 rational or integral polytope if all its vertices have rational or 
 integral coordinates, respectively. We write $d=\dim P$ and call $P$ 
 a $d$-polytope.

 For $P\subset {\mathbb  R}^m$ a rational $d$-polytope and $t$ a 
 nonnegative integer 
 we define the functions $i(P;t)=|tP\cap{\mathbb Z}^m |$ and
 $\overline{i}(P;t)=|tP^o\cap{\mathbb Z}^m |$, 
 where $tP=\{ta| a\in P\}$ is the $t$-fold dilatation of $P$. 
 Ehrhart proved in~\cite{Ehr67a, Ehr67b} that these functions are 
 \emph{quasi-polynomials},\index{polynomial!quasi-} 
 that is, they are of the form 
\[
    c_d(t)t^d+c_{t-1}(t)t^{d-1}+\ldots+c_0(t) ,
\]
 where  each $c_i(t)$ is a periodic function with integer period. Since 
 $i(P;t)$ is a quasi-polynomial, it can be defined for all $t\in
 {\mathbb Z}$. In fact, we have the following 
 reciprocity law due to Ehrhart~\cite{Ehr67c}:
\[
        i(P;-t)=(-1)^{d}\,\overline{i}(P;t) . 
\] 
 For more on this beautiful theory see the monograph by Ehrhart~\cite{Ehr77}.
 
 From Ehrhart~\cite{Ehr67a,Ehr67b} (also see Stanley~\cite{Sta96}),
 we have that when $P\subset {\mathbb  R}^m$ is an integral
 $d$-polytope,  
 then $i(P;t)$ and $\overline{i}(P;t)$ are polynomials, 
  which leads to the following definition of the Ehrhart polynomial.

\begin{definition} 
Let $P$ be an integral convex $d$-polytope.  Then the Ehrhart
polynomial of $P$ is  
  \[ 
       i(P;t)=c_0+c_1t+\ldots +c_{d-1}t^{n-1}+c_{d}t^n. 
 \] 
 \end{definition}    
 
From the early works of 
 Ehrhart~\cite{Ehr67a} and Macdonald~\cite{Mac71} it is known that
 $c_0=1$ and  
 %
 $c_d=\mbox{vol}(P)$,  and  that $c_{d-1}$ is half of the surface area 
 of $P$, normalized with respect to the sub-lattice on each face of 
 $P$.  Specifically, $c_{d-1}=1/2\sum_{F}\mbox{vol}_{d-1}(F)$, where $F$ 
 ranges over all facets of $P$ and the volume of a facet is measured 
 intrinsically with respect to the lattice ${\mathbb Z}^m\cap L_F$, 
 where $L_F$ is the affine hull of $F$.  The other coefficients were
 not well understood, until the later  
 work of Betke and Kneser~\cite{BK85}, Pommersheim~\cite{Pom93},
 Kantor and Khovanskii~\cite{KK93}  
 and Diaz and Robins~\cite{DR97}, but
 such interpretations go beyond  the scope of this chapter. For the 
 complexity of computing these coefficients see Barvinok~\cite{Bar94}.

 In the special case that $P$ is a \emph{zonotope} there 
 is a combinatorial interpretation for the coefficients of the Ehrhart
 polynomial.  First  
 recall that if  $A$ is an $r\times m$ real matrix 
 written in the form $A=[a_1,...,a_r]$,  then it defines a 
 \emph{zonotope}\index{zonotope} $Z(A)$ which  
 consists of those  points $p$ of ${\mathbb R}^m$ which can be 
 expressed in the form  
 \[ 
     p=\sum_{i=1}^m\lambda_ia_i, \quad 0\leq\lambda_i\leq 1. 
 \] 
 In other words, $Z(A)$ is the {\it Minkowski sum} of the line 
 segments  $[0,a_i]$, for $1\leq i\leq n$. For more on zonotopes, 
 see McMullen~\cite{McM71}.  
 
 When $A$ has integer entries, Stanley~\cite{Sta80}, using 
 techniques  from Shephard~\cite{She74}, proved that $i(P;t)=\sum_{X} 
 f(X)t^{|X|}$, where $X$ ranges over all linearly independent subsets 
 of columns of $A$ and where $f(X)$ denotes the greatest common 
 divisor of all minors of sizes $|X|$ of the matrix $A$. 
 
 When $A$ is a \emph{totally unimodular matrix},
\index{matrix!totally unimodular} that is, the determinant of every
square submatrix is 0 or $\pm 1$, then $Z(A)$ is described  
 as a \emph{unimodular zonotope}.\index{zonotope!unimodular}  For  
 these polytopes the previous result shows that 
 \[  
      i(Z(A);t)=\sum^r_{k=0}f_k\,t^k, 
 \] 
 where $f_k$ is the number of subsets of columns of the 
 matrix $A$ which are linearly independent and have cardinality $k$.  
 In other words, the Ehrhart polynomial $i(Z(A);t)$ is the generating 
 function of the number of independent sets in the regular matroid 
 $M(A)$. 
 
 The incidence matrix $D(G)$ of a graph $G$ is totally unimodular, a
 long-standing result  
 due to Poincar\'e~\cite{Poi01} with a modern treatment given by
 Biggs~\cite{Big96}. A 
 linearly independent subset of columns in $D$  
 corresponds to a subset of edges with no cycle. Thus, the coefficient
 $f_k$ in  
 this case is the number of spanning forests of $G$ with exactly $k$ 
 edges. From the previous chapter  we know that
 $T(G;x+1,y)=\sum_{k=0}^{r}f_{k}x^{r-k}$, where $r$ is the rank of the
 graph 
 $G$. With these ingredients we get the following relations with the 
 Tutte polynomial from Welsh~\cite{Wel97}.

 \begin{theorem} \label{reciprocity}
  If $G$ is a graph and $D$ is its incident matrix then the Ehrhart 
 polynomial of the unimodular zonotope $Z(D)$ is given by 
 \[ 
     i(Z(D);t)=t^{r}T(G;1+\frac{1}{t},1), 
 \] 
 where $r$ is the rank of $G$. 
 \end{theorem} 
 In this case, the zonotope $Z(D)$ is a 
 $r$-polytope in ${\mathbb  R}^n$, where $n$ is the number of vertices
 of $G$.   
 
 The reciprocity law of Theorem~\ref{reciprocity} leads to the
 following geometric result, 
 also from Welsh~\cite{Wel97}.
\begin{corollary} 
  If $D$ is the incidence matrix of a  rank $r$ graph $G$ with $n$ 
 vertices then for any positive integer $\lambda$ the number of  
 lattice points of ${\mathbb R}^n$ lying strictly inside the zonotope 
 $t Z(D)$ is given by  
 \[ 
     \overline{i}(Z(D);t)=(-t)^{r}T(G;1-\frac{1}{t},1). 
 \] 
\end{corollary} 
 In particular we have that the number of lattice points strictly
 inside $Z(D)$ is  $(-1)^{r}T(G;0,1)$.

\subsection{The Topological Tutte Polynomial of Bollob\'as and Riordan}
\index{polynomial!topological Tutte}  
 
  The classical Tutte polynomial discussed in the previous chapter is
  an invariant of abstract graphs, so it encodes  
  no information specific to graphs embedded in surfaces.  In
  \cite{BR01,BR02},  
  Bollob\'{a}s and Riordan generalize the classical Tutte polynomial to 
  topological graphs,\index{topological graphs}
  \index{graph!topological} that is, graphs embedded in
  surfaces. In \cite{BR01},  
  Bollob\'{a}s and Riordan define the {\em cyclic graph
  polynomial},\index{cyclic graph polynomial}
  \index{polynomial!cyclic graph} a three  
  variable deletion/contraction invariant for graphs embedded in 
  oriented surfaces.  They extend this work in \cite{BR02}, using a
  different approach,  
   with the four variable {\em ribbon graph polynomial}.
  \index{ribbon graph polynomial}\index{polynomial!ribbon graph}
  Both of these  polynomials extend the classical Tutte  
  polynomial, but in such a way that topological information about the 
  embedding is encoded.  The version for oriented surfaces is subsumed 
  by the version for arbitrary surfaces, so we focus on the latter
  here. The ribbon graph polynomial is also
  sometimes called the \emph{Bollob\'{a}s-Riordan polynomial} 
  \index{polynomial!Bollob\'{a}s-Riordan}after
  the authors or the \emph{topological Tutte polynomial} to emphasize
  that it simultaneously  encodes   topological information  
  while generalizing the classical Tutte polynomial.
   
\begin{figure}[hbtp] 
     \begin{center} 
         \includegraphics{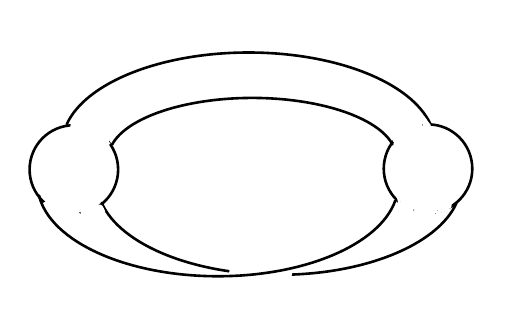} 
       \caption{A ribbon graph which is a M\"obius band with 
         boundary}\label{Fig:fatgraph}   
     \end{center} 
   \end{figure}

  First recall that a cellular embedding \index{cellular embedding}
  of a graph in an orientable or  
  unorientable surface can be specified by providing a sign for each 
  edge and a rotation scheme \index{rotation scheme}for the set of
  half edges at each vertex, where  
  a rotation scheme is simply a cyclic ordering of the half edges about a 
  vertex.  This is equivalent to a {\em ribbon} (or {\em fat}) {\em graph},
  \index{ribbon graph}\index{fat graph} which is a  
  surface with boundary where the vertices are represented by a set of 
  disks and the edges by ribbons, with the ribbon of an edge with a 
  negative sign having a half-twist. This can also be thought of as taking a 
  slight `fattening' of the edges of the graph as it is embedded in the 
  surface, or equivalently as `cutting out' the graph together with a small
  neighborhood of it  
  from the surface.  Figure~\ref{Fig:fatgraph} shows a graph with 
  two vertices and two  
  parallel edges, one positive and one negative.  It is embedded on a 
  Klein bottle, and the ribbon graph is a M\"obius band with boundary.  
 
  In addition to the usual graphic characteristics such as number of 
  vertices, connected components, rank, and nullity  for a ribbon
  graph $G$ we 
  also consider $bc\left( G \right)$, the number of boundary
  components \index{boundary components} of  
  the surface, and $t(G)$, an index of the
  orientability of the  
  surface.  The value of $t(G)$ is 0 if the surface is
  orientable, and 1  
  if it is not.  Thus, $t(G)$ is 1 if and only if  for
  some cycle in $G$,  
  the product of the signs of the edges is negative.  
   
   
   \begin{definition}\label{ribbon graph poly}
  \index{topological Tutte polynomial}
  \index{polynomial!topological Tutte} Let $G$ be a ribbon
  graph, that is, a graph  
  embedded in a surface. Then the topological Tutte polynomial of
  Bollob\'as and Riordan is given by 
  \begin{equation*} 
   \begin{split} 
   R(G;x,y,z,w) &= \sum_{A \subseteq E( G)} 
                 (x - 1)^{r( G ) - r( A )} 
                 y^{n(A)} z^{\kappa(A) - bc(A) + n(A)} 
                   w^{t(A)}\\ 
                &\in Z[x,y,z,w]/\vphantom {Z[x,y,z,w]} 
                  \langle w^2  - w\rangle . 
   \end{split} 
  \end{equation*} 
 
  \end{definition} 
 
  As previously, $r(A)$, $\kappa(A)$, $n(A)$, and now also
  $bc(A)$ and  
  $t(A)$, refer to the spanning subgraph of $G$ with
  edge set $A$,   
  here with its embedding inherited from $G$.  
   
  Clearly, by comparing with the rank-nullity generating function
  definition of the classical Tutte polynomial given in the previous
  chapter, this generalizes the  
  classical Tutte polynomial.  Like the classical Tutte polynomial, 
  $R(G;x,y,z,w)$ 
  is multiplicative on disjoint 
  unions and one point joins of ribbon graphs. More importantly, it 
  retains the essential properly of obeying a deletion/contraction
  reduction relation.  
   
  We must first define deletion and contraction in the context of
  embedded graphs. The ribbon graph resulting from deleting an edge is  
  clear, but contraction requires some care.  Let $e$ be a non-loop 
  edge.  First assume the sign of $e$ is positive, by flipping one 
  endpoint if necessary to remove the half twist (this reverses 
  the cyclic order of the half edges at that vertex and toggles their
  signs). Then $G/e$ is formed by  
  deleting $e$ and identifying its endpoints into a single vertex $v$. 
  The cyclic order of edges at $v$ comes from the original cyclic order 
  at one endpoint, beginning where $e$ had been, and continuing with the 
  cyclic order at the other endpoint, again beginning where $e$ had 
  been.   
 
   
  \begin{theorem}\label{ribbon delete-contract} 
  If $G$ is a ribbon graph, then  
   \[ R(G;x,y,z,w) = R(G/e;x,y,z,w) + R(G-e;x,y,z,w)\]
   if $e$ is an ordinary edge and $R(G;x,y,z,w) =$
  $xR(G/e;x,y,z,w)$ if   $e$ is a bridge.  
  \end{theorem} 
   
  The proof depends on a careful analysis of how each of the relevant
  parameters  
  $r(A)$, $\kappa(A)$, $n(A)$, $bc(A)$ and
  $t(A)$ changes with the deletion  
  or contraction of an edge.  
   
  Repeated application of this theorem reduces a ribbon graph to a 
  disjoint union of embedded \emph{blossom graphs},
  \index{blossom graph}\index{graph!blossom} that is, graphs each  
  consisting of a single vertex with some number of loops.  Because of 
  the embedding, the loops are signed, and there is a rotation system of 
  half-edges about the single vertex.  Not surprisingly, the topological 
  information is distilled into these minors of the original graph, and 
  to complete a deletion/contraction linear recursion computation, it
  is necessary to specify an evaluation of these  
  terminal forms.  
   
\begin{figure}[hbtp] 
     \begin{center} 
       \includegraphics{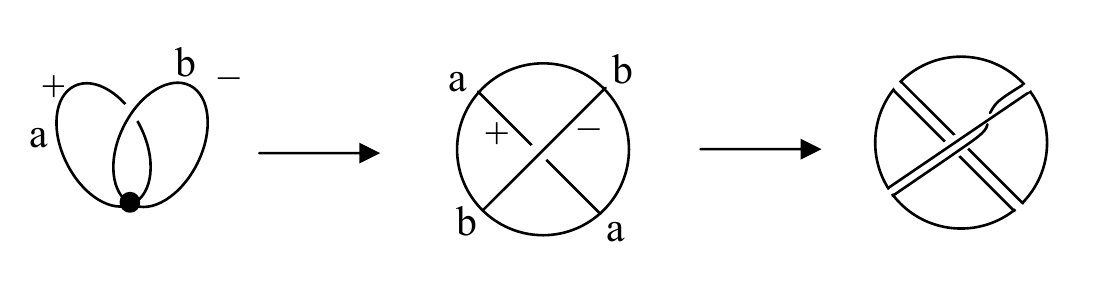} 
       \caption{A signed blossom graph, its signed chord diagram, and
         the boundary components of the signed chord
         diagram}\label{Fig:figureB}    
     \end{center} 
   \end{figure}

  Signed chord diagrams \index{chord diagram} provide a useful device
  for determining the  
  relevant parameters of an embedded blossom graph.  Recall that a
  \emph{chord  diagram} consists of a circle with $n$ symbols on its
  perimeter, with  
  each symbol appearing twice and a chord drawn between each pair of 
  like symbols.  A \emph{signed chord diagram}
  \index{chord diagram!signed}\index{signed chord diagram} simply has
  a sign on   each chord.  
  A signed chord diagram $D$ corresponds to an embedded blossom graph 
  $G$ by assigning a symbol to each loop and arranging them on the 
  perimeter of the circle in the chord diagram in the same order as
  the cyclic order of the half-edges  
  about the vertex.  A chord receives the same sign as the loop it 
  represents.  If we `fatten' the chords as in 
  Figure~\ref{Fig:figureB}, with a negative  
  chord receiving a half-twist, then $bc(G)$ is equal to the number of 
  components in the resulting diagram, which is denoted
  $bc(D)$.  
  Similarly, since $G$ has only one vertex, $n(G)$ is the number of 
  edges of $G$, which is the number of chords of $D$, so we denote this 
  by $n(D)$.  We also set
  $t(D)=t(G)$, and note that
  $t (D) = 0$   if all  
  chords of $D$ have a positive sign, and $t(D)=1$
  otherwise.  This,  
  combined with the definition of $R(G;x,y,z,w)$ above, gives the 
  following evaluation for these terminal forms.  

  \begin{theorem}\label{blossom eval} 
   If $G$ is an embedded blossom graph with 
  corresponding signed chord diagram $D$, then  
  \[ 
   R(G;x,y,z,w) = \sum_{D' \subseteq D} y^{n(D')}   
                                      z^{1- bc(D')+ n(D')} 
                                      w^{t(D')}, 
  \] 
  where the sum is over all subdiagrams $D'$ of $D$.  
  \end{theorem} 
%
%
%
   
Theorems~\ref{ribbon delete-contract} and~\ref{blossom eval} taken
together give a linear recursion definition for $R(G;$ $x$, $y$, $z$, $w)$.  
There are a number of technical considerations, similar to the care 
  that must be taken in contracting edges, but nevertheless many other 
  properties analogous to those of the classical Tutte polynomial hold. For 
  example, $R(G;x,y,z,w)$ has 
  a spanning tree expansion, a universality property, and duality 
  relation (in addition to the Bollob\'{a}s and Riordan's originating work 
  in \cite{BR01, BR02}, see also recent work by Chmutov \cite{Chm} and
  Moffatt  
  \cite{Mof08}).  Furthermore, Chmutov and Pak \cite{CP07}, and 
  Moffatt \cite{Mof,Mof08} have shown that $R(G;x,y,z,w)$ also extends   
  the relation between the classical Tutte polynomial and the Kauffman 
  bracket and Jones polynomial of knot theory due to
  Thistlethwaite~\cite{Thi87} and Kauffman~\cite{Kau89}.

 \subsection{Martin, or Circuit Partition,
   Polynomials}\label{sec:cpp} 
 
In his 1977 thesis, Martin~\cite{Mar77} recursively defined
polynomials $M(G,x)$ and $m(\vec{G};x)$ that encode, respectively,
information about the families of circuits in 4-regular Eulerian
graphs and digraphs. Las Vergnas subsequently found a state model
expression for these polynomials, extended their properties to general
Eulerian graphs and digraphs, and  
further developed their theory (see \cite{Las79, Las83, Las88}).  Both
Martin~\cite{Mar78} and Las Vergnas~\cite{Las88} found combinatorial
interpretations for some small integer evaluations of the polynomials,
while combinatorial interpretations for all integer values as well as
some derivatives were given in \cite{E-M00, E-M04a, E-M04b}, and by
Bollob\'{a}s~\cite{Bol02}.    
 
Transforms of the Martin polynomials, $J(G;x)$ and $j(\vec{G};x)$,
given in \cite{E-M98}, and then aptly named  
circuit partition polynomials in \cite{ABS00}, facilitate these
computations, and for this reason we give the definitions below in
terms of $J$ and $j$.   
Like many of the polynomials surveyed here, the circuit partition
polynomials have several definitions, including linear recursion
formulations, generating function formulations, and state model
formulations.  We give the state model definition, and refer the
reader to ~\cite{E-M04a, E-M04b} for the others.  
 
As with other state model formulations, we must first specify what we
mean by a state of a graph (or digraph) in this context.  Here an
Eulerian graph \index{graph!Eulerian} must have vertices all of even
degree, but it need not be connected.  An Eulerian digraph
\index{Eulerian digraph}\index{graph!Eulerian digraph} must have the
indegree equal to the outdegree at each vertex, and again need not be
connected.  
 
\begin{definition} \label{cpp graph state}
 An \textit{Eulerian graph state}\index{Eulerian graph state}
 of an Eulerian graph $G$ is the result of replacing each
$2n$-valent vertex $v$ of $G$ with $n$  
2-valent vertices joining pairs of half edges originally adjacent to $v$. 
An Eulerian graph state of an Eulerian digraph $\vec {G}$ is defined 
similarly, except here each incoming half edge must be paired with an
outgoing  edge. 
\end{definition} 

Note that a Eulerian graph state is a disjoint union of cycles, each
consistently oriented in the case of a digraph.  
 
\begin{definition} The circuit partition polynomial. 
   \index{circuit partition polynomial}
   \index{polynomial!circuit partition} Let $G$
  be an Eulerian graph, let $St(G)$ be the set of states of $G$, and
  let $c(S)$ be the number of components in a state $S \in St(G)$.
  Then the circuit partition polynomial has a state model formulation
  given by   
\[ 
J(G;x) = 
\sum\limits_{S \in St(G)}{x^{c(S)}}.  
\] 

The circuit partition polynomial is defined similarly for Eulerian
digraphs as   
\[ 
j(\vec {G};x) = \sum\limits_{S \in St(\vec{G})}{x^{c(S)}}.  
\] 
\end{definition}

The transforms between the circuit partition polynomials and the original
Martin polynomial, as extended to general Eulerian graphs and digraphs
by Las Vergnas, are: 
 
\begin{equation}
 J(G;x) = xM\left( {G;x + 2} \right), \text{ for }G \text{ an Eulerian
 graph, and }   
\end{equation} 
\begin{equation}\label{jm} 
j(\vec {G};x) = xm\left( 
{\vec {G};x + 1} \right) \text{ for } \vec {G} \text{ an Eulerian digraph}. 
\end{equation} 
 
The circuit partition polynomials have `splitting' formulas, analogous
to Tutte's identity for the chromatic polynomial given in the previous chapter, proofs for
which may be found in \cite{E-M98, E-M04b}.  These formulas derive
from the Hopf algebra structures of the generalized transition
polynomial discussed in Subsection~\ref{gen_transition_poly},
but may also be proved combinatorially, as by
Bollob\'as~\cite{Bol02} and \cite{E-M04b}.  
 
\begin{theorem}\label{cppsplit} Let $G$ be an Eulerian graph and
  $\vec{G}$ be an Eulerian digraph.  Then  
 
\[ 
J(G;x + y) = \sum {J\left( {A;x} 
\right)J\left( {A^c;y} \right)}, 
\] 
where the sum is over all subsets $A \subseteq E(G)$ such that $G$
restricted to both $A$ and $A^c = E(G)-A$ is Eulerian. Also,   
\[ 
j(\vec {G};x + y) 
= \sum {j\left( {\vec {A};x} \right)j\left( {\vec {A}^c;y} \right)}, 
\] 
where the sum is over all subsets $\vec {A} \subseteq E(\vec {G})$
such that $\vec {G}$ restricted to both $\vec {A}$ and $\vec {A}^c$ is
an Eulerian digraph.  
\end{theorem}

The connection between the circuit partition polynomial of a digraph
and the Tutte polynomial of a planar graph $G$ is through the oriented
medial graph $\vec{G}_{m}$ described in the previous chapter.
Martin~\cite{Mar77} proved the following, which we extend to the
circuit partition polynomial via (\ref{jm}).  
 
\begin{theorem}\label{Tutte-Martin} 
  
Let $G$ be a connected planar graph, and let $\vec {G}_m$ be its oriented 
medial graph. Then relationships among the Martin polynomial, circuit
partition polynomial, and Tutte polynomial are: 
\[
j(\vec {G};x)= x m(\vec {G}_m ;x+1) = xt(G;x+1,x+1). 
\] 
\end{theorem} 
 
\begin{figure}[hbtp] 
     \begin{center} 
         \includegraphics[width=0.9\textwidth]{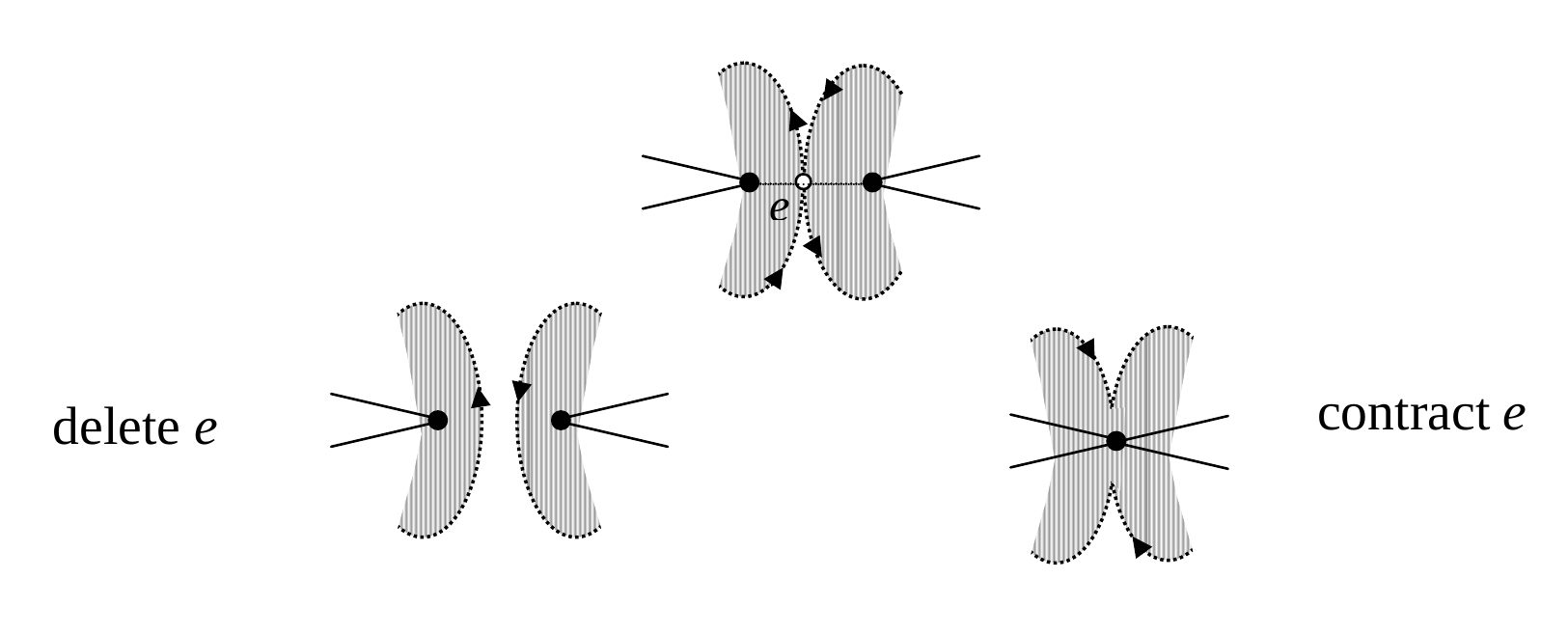} 
       \caption{An edge $e$ in a planar graph $G$, with the
         corresponding vertex $v$ in the oriented medial graph
         $\vec{G}_m$ (dotted edges).  Deleting $e$ corresponds to one
         possible configuration at $v$ in an Eulerian graph state of
         $\vec{G}_m$, while contracting $e$ corresponds to the
         other}\label{Fig:TutteMartinfig}     
     \end{center} 
   \end{figure}

The proof of this theorem depends on a fundamental observation
relating deletion/contraction in $G$ with choices of configurations at
a vertex in an Eulerian graph state of $\vec {G}_m $, as illustrated
in Figure~\ref{Fig:TutteMartinfig}.  Theorems \ref{cppsplit} and
\ref{Tutte-Martin} combine to give the basis for many of the
combinatorial interpretation of the Tutte polynomial along the line $y
= x$ described in the previous chapter.  For more details, see 
Martin~\cite{Mar77, Mar78}, Las Vergnas~\cite{Las79, Las83, Las88},
Bollob\'as~\cite{Bol02}, and also~\cite{E-M98, E-M00,E-M04a,
  E-M04b}.
 
Evolving from the relation between the Tutte and Martin polynomials is
the theory of isotropic systems, \index{isotropic system} which
unifies essential properties of 4-regular graphs and pairs of dual
binary matroids. A series of papers throughout the 1980's and 1990's,
including work by Bouchet \cite{Bou87a}, \cite{Bou87b}, \cite{Bou88},
\cite{Bou89}, \cite{Bou91},\cite{Bou93}, as well as Bouchet and Ghier
\cite{BG96}, and Jackson \cite{Jac91}, significantly extends the
relationship between the Tutte polynomial of a planar graph and the
Martin polynomial of its medial graph via the theory of isotropic
systems.  
 
 \subsection{Interlace Polynomial} 
 
In \cite{ABS00}, Arratia, Bollob\'as and Sorkin 
defined a one-variable graph polynomial motivated by questions 
arising from DNA sequencing by hybridization 
addressed by Arratia, Bollob\'as, Coppersmith and Sorkin in 
\cite{ABCS00}, an application we will return to in
Section~\ref{sec:applications}.   
 In \cite{ABS04b}, Arratia, Bollob\'as, and Sorkin defined a 
two-variable interlace polynomial, \index{interlace polynomial} and
showed that the original polynomial of  
\cite{ABS00} is a specialization of it, renaming the original one-variable 
polynomial as the vertex-nullity interlace polynomial
\index{vertex-nullity interlace polynomial}
\index{polynomial!vertex-nullity} due to its  
relationship with the two-variable generalization.  
 
Remarkably, despite very different terminologies, motivations and
approaches, the original vertex-nullity interlace polynomial of a
graph may be realized as the Tutte-Martin polynomial
\index{Tutte-Martin polynomial}\index{polynomial!Tutte-Martin} of an
associated isotropic system (see Bouchet
\cite{Bou05}). For exploration of this relationship, see the works
mentioned in Subsection \ref{sec:cpp}, as well as Aigner \cite{Aig00},
Aigner and Mielke \cite{AM00}, Aigner and van der Holst \cite{AvdH04},
Allys \cite{All94}, and also Bouchet's series on multimatroids
\cite{Bou97,Bou98a,Bou98b,Bou01}.   
 
\begin{figure}[hbtp] 
     \begin{center} 
         \includegraphics[width=0.9\textwidth]{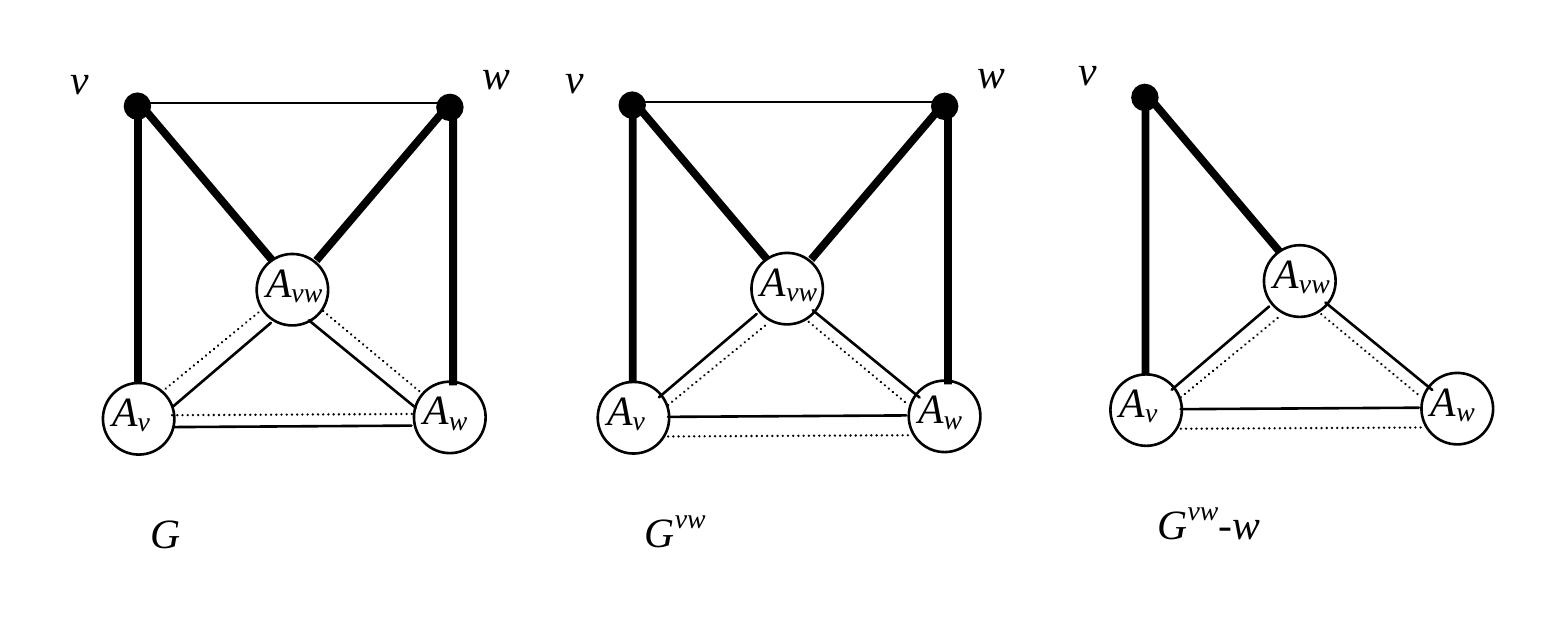} 
       \caption{Pivoting on the edge $vw$. $A_{v}$,  $A_{w}$ and
         $A_{vw}$  are the  sets of vertices of $G$ 
         adjacent to $v$ only, $w$ only, and to both $v$ and $w$,
         respectively.  These sets are constant in all the diagrams.
         Vertices of $G$ adjacent to neither $v$ nor $w$ are omitted.  Heavy
         lines indicate that all edges are present, and dotted lines
         represent non-edges. Note interchange of edges and non-edges
         among $A_{v}$, $A_{w}$ and $A_{vw}$}\label{Fig:pivotfig}     
     \end{center} 
   \end{figure}

Both the vertex-nullity interlace polynomial of a graph and the
two-variable interlace  
polynomial may be defined recursively via a  
pivot operation. \index{pivot operation} This pivot is defined as
follows.  Let $vw$ be an edge of a graph $G$, and let $A_v$,  
$A_w$ and $A_{vw}$ be the sets of vertices 
in $V(G)\setminus \{v,w\}$ 
adjacent to $v$ only, $w$ only, 
and to both $v$ and $w$, respectively.  The pivot operation ``toggles" 
the edges among $A_v$,  $A_w$ and $A_{vw}$, by deleting existing
edges and inserting  
edges between previously non-adjacent vertices. 
The result of this operation 
is denoted $G^{vw}$.  More formally,  $G^{vw}$  has the same vertex
set as $G$, and edge  
set equal to the symmetric difference $E(G)\Delta   S$, where $S$  is the 
complete tripartite graph with vertex classes 
$A_v$,  $A_w$ and $A_{vw}$. See Figure~\ref{Fig:pivotfig}.\\ 
 
Also, $G^a$ is the {\em local complementation}
\index{local complementations} of $G$, defined   
as follows. Let $N(a)$ be the neighbors of $a$, that is, the set 
$\{w\in V: a {\text{ and }}w {\text{ are joined by an edge}}\}$. 
The graph $G^a$ is equal to $G$ except 
that  
we ``toggle" the edges among the neighbors of $a$, 
switching edges to non-edges and vice-versa. 
 
\begin{definition} 
Let $G$ be a graph of order $n$, which may have loops, but no multiple
loops or multiple edges.  The two-variable interlace polynomial may be
given recursively by $q(E_n)=y^n$ for $E_n$, the edgeless graph  on
$n\geq 0$ vertices, with   
\[
q(G)=q(G- a)+q(G^{ab}- b)+((x-1)^2-1)q(G^{ab}- a- b), \\ 
\]  
\noindent 
for any edge $ab$ where neither $a$ nor $b$ has 
a loop, and   
\[
q(G)=q(G- a)+(x-1)q(G^a- a), \\ 
\]
\noindent 
for any looped vertex $a$.  
\end{definition} 
 
Alternatively, the interlace polynomial has the following generating
function representation.  
 
\begin{definition} Let $G$ be a graph of order $n$, which may have
  loops, but no multiple loops or multiple edges.  Then the
  two-variable interlace polynomial may be given by  
 
\begin{equation*}
  q(G;x,y)=\sum_{S\subseteq
  V(G)}(x-1)^{r(\left. G \right|_S)}(y-1)^{n(\left. G
  \right|_S)},\end{equation*}  
\noindent  
where $r(\left. G \right|_S)$ and $n(\left. G
\right|_S)=|S|-r(\left. G \right|_S)$ are, respectively,  
the ${\mathbb F}_2$-rank and 
nullity of the adjacency 
matrix of $\left. G \right|_S$, the subgraph of $G$ restricted to $S$. 
\end{definition}

\begin{definition}\label{vertexnullity} 
The vertex-nullity interlace polynomial is defined recursively as: 
$$ 
q_N(G;x)= 
\begin{cases} 
x^n {\text{ if }}G=E_n, {\text{ the edgeless graph on $n$ vertices}}\\ 
q_N(G- v;x)+q_N(G^{vw}- w;x) {\text{ if }}vw\in E(G). 
\end{cases} 
$$ 
 
\end{definition} 
 
This polynomial was shown to be well defined by Arratia, Bollob\'as,
and Sorkin for all simple graphs in \cite{ABS00}, and then was shown
in \cite{ABS04b} to be a specialization of the two-variable interlace
polynomial as follows.  
 
\begin{equation*}
q_N(G;y)=q(G;2,y)=\sum_{W\subseteq V(G)}(y-1)^{n(\left. G \right|_W)}. 
\end{equation*}

An equivalent formulation for $q_N(G;x)$  for simple graphs is given
by Aigner and van der Holst in  
\cite{AvdH04}. 
 
A somewhat circuitous route through the circuit partition polynomial
relates the vertex-nullity interlace polynomial to the Tutte
polynomial. \index{Tutte polynomial}\index{polynomial!Tutte} First
recall that a {\em circle graph}\index{circle graph} on $n$ vertices
is a graph $G$  
derived from a chord diagram. 
Two vertices $v$ and $w$ in $G$ share an edge if and only if their
corresponding  
chords intersect in the chord diagram. Note that $G$ is necessarily simple. 
 
For circle graphs, the vertex-nullity 
interlace polynomial and the circuit partition polynomial 
are related 
by the following theorem, noting that  
although $\vec{G}$ may be a multigraph, $H$ is necessarily  
simple. 
 
\begin{theorem}(Arratia, Bollob\'as and Sorkin \cite{ABS00}, Theorem
  6.1).\label{A}
If $\vec{G}$  is a $4$-regular Eulerian digraph, $C$ is any Eulerian 
circuit of $\vec{G}$, and  $H$ is the circle graph of the chord diagram 
determined by $C$, 
 then $j(\vec{G};x)=xq_N(H;x+1)$. 
\end{theorem} 
 
This now allows 
us to relate the vertex-nullity interlace polynomial 
to the Tutte polynomial, \index{Tutte polynomial}\index{polynomial!Tutte}
a relation proved in \cite{E-MS07} and also observed by Arratia,
Bollob\'as and Sorkin at the end of  
Section $7$ in \cite{ABS04b}.

\begin{theorem}\label{B} 
If $G$ is a planar graph, and $H$ is the circle graph of some Eulerian 
circuit of  $\vec{G}_m$, 
then $q_N(H;x)=t(G;x,x)$. 
\end{theorem} 
 
\begin{proof} 
By Theorem~\ref{A}, $j(\vec{G}_m;x)=xq_N(H;x+1)$, but 
recalling that the circuit partition and Martin polynomials are 
simple translations of each other, 
we have from Theorem~\ref{Tutte-Martin} that
$j(\vec{G}_m;x)=xm(\vec{G}_m;x+1)$, and hence  
$q_N(H;x)=m(\vec{G}_m;x)=t(G;x,x)$.  
\end{proof} 
 
The interlace polynomial has generated further interest and other
applications in  
Balister, Bollob\'as, Cutler, and Pebody \cite{BBCP02}, 
and Balister, Bollob\'as, Riordan, and Scott \cite{BBRS01}, Glantz
and Pelillo \cite{GP06}, Ellis-Monaghan and Sarmiento \cite{E-MS07}.

\section{Multivariable Extensions}
\label{sec:multivariable}

 Multivariable extensions have proved valuable theoretical tools for
 many of the polynomials we have seen since they capture information
 not encoded by the original polynomial.  
More critically, powerful algebraic tools not applicable to the
original polynomial may be available to the multivariable version,
providing new means of extracting combinatorial information from the
polynomial.  While the multivariable indexing may make the defining
notation somewhat bulky, these generalizations are natural extensions
of classical versions, computed in exactly the same ways, only now
also keeping track of some additional parameters in the computation
processes.

\subsection{Generalized coloring polynomials \\and the U-polynomial} 
\label{subsec:symmetric}

  The evaluation of the chromatic polynomial at $\lambda$ can be written as 
  \begin{equation}\label{chrompoly}
  \index{chromatic polynomial}\index{polynomial!chromatic}  
  \chi_{G}(\lambda)=
      \sum_{\substack{\phi:V\rightarrow \{1,\ldots,\lambda\}\\ 
      \text{proper}}}{1}. 
  \end{equation} 
   
  This was generalized to a symmetric function  over (commuting) 
  indeterminates  $x_1,x_2,\ldots$ by Stanley~\cite{Sta95} in the
  following way.   
 
\begin{definition}\label{symmetric chromatic} 
Let $G=(V,G)$ be a graph, let $\phi:V\rightarrow
   \mathbb{P}=\{1,2,\ldots\}$, and denote  
   the product $\prod_{v\in V}x_{\phi(v)}$ by $x^{\phi}$. Then the
   symmetric function generalization of the chromatic polynomial is
   \index{symmetric chromatic polynomial}
   \index{polynomial!symmetric chromatic}
  \[ 
  X_{G}(\boldsymbol{x})= X(G;x_1,x_2,\ldots)= 
  \sum_{\substack{\phi:V\rightarrow  \mathbb{P}\\ \text{proper}}} 
  {x^{\phi}}. 
  \] 
  \end{definition} 
 
  That this is a generalization of the chromatic polynomial can be
  seen by setting $x_i=1$ for $1\leq i\leq \lambda$ and  $x_j=0$ for
  $j> \lambda$ and noting that the expression in (\ref{chrompoly}) for
  the chromatic polynomial evaluated at $\lambda$ results.    
   
  Generalizing polynomial graph invariants is not a theoretical
  exercise. The original invariant encodes combinatorial information,
  and the multivariable generalization will encode not only the same
  information but also more refined information.  For example, the
  chromatic polynomial  
  of any tree with $n$ vertices has  
  chromatic polynomial $x(x-1)^{n-1}$. But not all trees  have the same 
  $X_{G}(\boldsymbol{x})$. For example, if $K_{1,3}$ is the 4-star 
  graph and $P^{4}$ is  
  the path of order 4,  $X_{K_{1,3}}(\boldsymbol{x})$ 
  has a term $x_{i}x_{j}^3$ for all $i\neq j$, but such a term is not 
  present in $X_{P^{4}}(\boldsymbol{x})$ . In fact, it is still an
  open question if $X$  
  distinguishes trees, that is if  $X_{T_1}(\boldsymbol{x})\neq 
  X_{T_2}(\boldsymbol{x})$, whenever $T_1$ and $T_2$ are not
  isomorphic trees.  
   
  A similar multivariable extension of the bad coloring polynomial is
  also natural, especially given the importance of the latter because
  of its being  
  equivalent to the Tutte polynomial. \index{Tutte polynomial}
  \index{polynomial!Tutte}  The following  
  generalization of the bad coloring polynomial is also due to 
  Stanley~\cite{Sta98}.
 
\begin{definition}\label{multivariable bad coloring}
\index{bad coloring polynomial!symmetric function extension}
\index{polynomial!symmetric bad coloring}  
Let $G=(V,E)$ be a graph, let $\phi:V\rightarrow
  \mathbb{P}=\{1,2,\ldots\}$, and let $b(\phi)$ be the set of
  monochromatic 
  edges in the coloring given by $\phi$. Then the symmetric function
  generalization of the bad coloring polynomial   
  over indeterminates $x_1,x_2,\ldots$ and $t$ is   
\[ X_{G}(\boldsymbol{x},t)= \sum_{\phi:V\rightarrow 
  \mathbb{P}}{(1+t)^{|b(\phi)|}x^{\phi}},
\] 
where the sum is over all possible colorings $\phi$ of the graph $G$.
\end{definition} 
 
  Again, by setting  $x_i=1$ for $1\leq 
  i\leq \lambda$ and  $x_j=0$ for $j> \lambda$ we get the bad-coloring 
  polynomial, and hence the Tutte polynomial. \index{Tutte polynomial}
  \index{polynomial!Tutte} Therefore,  
  $X_{G}(\boldsymbol{x},t)$ is a multivariable generalization of the Tutte 
  polynomial.  
 
 There is another multivariable generalization of the Tutte polynomial
  that was developed  
  independently and for very different reasons.   This generalization is  
  called the U-polynomial and is due to Noble and Welsh in~\cite{NW99}.   
 
\begin{definition} \label{U-poly}\index{U-polynomial}
\index{polynomial!U-}
Let $G=(V,E)$ be a graph.  Then the U-polynomial of $G$ is 
\[  
   U_{G}(\boldsymbol{x},y)= \sum_{A\subseteq E}{x_{n_1}\cdots x_{n_k} 
   (y-1)^{|A|-r(A)}}, 
\] 
  where $n_1,\ldots,n_k$ are the numbers of vertices in the $k$ different 
  components of $G$ restricted to $A$.  
\end{definition}

  Clearly, this is a generalization of the Tutte 
  polynomial, as by setting $x_i=(x-1)$ for all $i$ in 
  $U_{G}(\boldsymbol{x},y)$ we get
  $(x-1)^{\kappa(G)}T_{G}(x,y)$. Note that the factor $x_{n_1}\cdots
  x_{n_k}$ in every term keeps track of the number of 
  vertices in the different components in $A$. This is a refinement of
  the rank-nullity generating-function definition of the Tutte polynomial
  where the factors $x^{r(G)-r(A)}=x^{\kappa(A)-\kappa(G)}$ in each
  term keep track of the number of components in $A$.

  That $U_G$ captures more combinatorial information from $G$ than the
  Tutte polynomial can be
  seen by noting that $U_G$ contains the matching generating
  polynomial,\index{matching generating polynomial} and thus the
  matching polynomial,\index{matching polynomial}
  as  a specialization as well.  
  \begin{theorem} 
    For any graph $G$,  
    \[ g(G;x)= U_{G}(1,t,0,\ldots,0,\ldots,y=1).\] 
  \end{theorem} 

 The U-polynomial has a deletion/contraction  reduction relationship 
 not in the class of graphs but in the class of weighted graphs. To
 see this, we turn to the W-polynomial\index{W-polynomial}
 \index{polynomial!W-} also due to Noble and Welsh
 in~\cite{NW99}.  A 
 \emph{weighted graph}\index{graph!weighted} consists of a graph
 $G=(V,E)$, together with a  weight function  $\omega:V\rightarrow
 \mathbb{Z}^{+}$. 

 If $e$ is an edge of $(G,\omega)$ then  $(G\setminus e,\omega)$ is the
 weighted graph obtained from $(G,\omega)$ by deleting $e$ and leaving
 $\omega$ unchanged. If $e$ is not a loop, $(G/e,\omega/e)$ is the
 weighted graph obtained from $(G,\omega)$ by  contracting $e$, that
 is deleting $e$ and identifying its endpoints $v$, $´v'$ into a
 single vertex $v''$. The weight function $\omega/e$ is defined as
 $\omega/e(u)=\omega(u)$ for all $u\in V\setminus \{v,v'\}$ and
 $\omega/e(v'')=\omega(v)+\omega(v')$.

\begin{definition} Let $(G,\,\omega)$ be a weighted graph.  The
 W-polynomial may be given recursively by the following rules. If
 $e$ is an ordinary edge or a bridge, then 
\[   W_{(G,\, \omega)}(\boldsymbol{x},y)= 
              W_{(G\setminus e,\, \omega)}(\boldsymbol{x},y) +
              W_{(G/e,\, \omega/e)}(\boldsymbol{x},y) .
\]
 If $e$ is a loop, then $ W_{(G,\,\omega)}(\boldsymbol{x},y)=yW_{(G\setminus
 e,\, \omega)} (\boldsymbol{x},y)$. Otherwise, $(G,\,\omega)$  is $E_n$, the
 edgeless graph  on  $n\geq 0$ vertices, with weights $a_{1},\ldots
 ,a_{n}$ and $W_{(E_n,\,\omega)}(\boldsymbol{x},y)=x_{a_1}\cdots x_{a_n}$.
\end{definition}

 That the resulting multivariate polynomial  W is independent of the
 order in which the edges are deleted and contracted is proved
 in~\cite{NW99}. This can easily  be done by induction on the number of edges
 once it is proved that  the order in which you contract or delete
 edges in $(G,\omega)$ does  not affect the weighted graph which you
 obtain.  

 The U-polynomial is obtained from the 
 W-polynomial by setting all weights equal to 1 and a proof that
 this definition is equivalent to Definition~\ref{U-poly} can be
 found in~\cite{NW99}. Actually in~\cite{NW99} it is proved that $W$ has a
 representation of the form
 \[ 
   W_{(G,\,\omega)}(\boldsymbol{x},y)=
             \sum_{A\subseteq E}{x_{c_1}\cdots  x_{c_k}  
                                      (y-1)^{|A|-r(A)}}  ,
 \]
where $c_i$, $1\leq i\leq k$, is the total weight of the $i$th
component of the weighted subgraph $(A,\omega)$.

 Noble and Welsh \cite{NW99} show that the symmetric function
 generalization of the bad coloring polynomial and the U-polynomial are
 equivalent in the following  sense.   
 
\begin{theorem} 
 For any graph $G$, the polynomials $U_G$ and $X_G$ determine each
 other in that if $p_0=1$ and $p_r=\sum_{i}{x_{i}^{r}}$, then    
\[ X_{G}(\boldsymbol{x},t)= t^{|V|}U_{G}(x_j=\frac{p_j}{t},y=t+1).\] 
\end{theorem}

 There is yet another polynomial, the \emph{polychromate}, introduced originally
 by Brylawski in~\cite{Bry81}, that is as general as $U_{G}$ or
 $X_{G}$. Given a graph $G$ and a partition $\pi$ of its vertices into
 non-empty blocks, let $e(\pi)$ be the number of edges with both ends
 in the same block of the partition. If $\tau(\pi)=(n_1,\ldots, n_k)$
 is the type of partition $\pi$, we denote by
 $\boldsymbol{x}_{\tau(\pi)}$  the monomial
 $\prod_{i=1}^{k}{x_i^{n_i}}$.

\begin{definition}

Let $G$ be a graph. Then the polychromate
 \index{polychromate} $\chi_{G}(\boldsymbol{x},y)$   is 

\[
          \chi_{G}(\boldsymbol{x},y)=\sum_{\pi}{y^{e(\pi)}
                                     \boldsymbol{x}_{\tau(\pi)}} ,
\]
where the sum is over all partitions of $V(G)$.  

\end{definition}

We have the following
theorem due to Sarmiento in~\cite{Sar00} but see~\cite{MN} for a
different proof. 
\begin{theorem} The polynomials $U_{G}(\boldsymbol{x},y)$ and
  $\chi_{G}(\boldsymbol{x},y)$ are equivalent.
\end{theorem}

 The story doesn't end here. All three polynomials
 $U_{G}(\boldsymbol{x},y)$, 
 $X_{G}(\boldsymbol{x},t)$ and $\chi_{G}(\boldsymbol{x},y)$ have
 natural extensions. For example the extension of the
 $X_{G}(\boldsymbol{x},t)$ replaces the $t$ variable by countably
 infinitely many variables $t_1,t_2,\ldots$, enumerating not just the
 total number of monochromatic edges but the number of monochromatic
 edges of each color. It is defined as follows.
 \[ X_{G}(\boldsymbol{x},\boldsymbol{t})= \sum_{\phi:V\rightarrow \mathbb{P}}
           { \left(\prod_{i=1}^{\infty}{(1+t_i)^{|b_i(\phi)|}}\right)
             x^{\phi}
           } ,
\] 
 where the sum is over all colorings $\phi$ of $G$ and $b_i(\phi)$ is
 the set of monochromatic edges for which both end points have
 color $i$. By setting $t_i=t$ for all $i\geq 1$ we regain
 $X_{G}(\boldsymbol{x},t)$.

 For the other extension the reader is referred to~\cite{MN}. There
 it is also proved that all these extensions are equivalent. 
  
  \subsection{The Parametrized Tutte Polynomial}\label{subsec:param tutte} 
 
  The basic idea of a \emph{parametrized Tutte polynomial} 
 \index{parametrized Tutte polynomial}
 \index{polynomial!parametrized Tutte} is to allow each  
  edge of a graph to have four parameters (four ring values specific to 
  that edge), which apply as the Tutte polynomial is computed via a 
  deletion/contraction recursion.  Which parameter is applied in a
  linear recursion reduction depends on  
  whether the edge is deleted or contracted as an ordinary edge, or 
  whether it is contracted as an isthmus or deleted as a loop.  
  The difficulty lies in assuring that a well defined function, that is, 
  one independent of the order of deletion/contraction, results.  This 
  requires a set of relations, coming from three very small graphs, to 
  be satisfied.  Interestingly, additional constraints are necessary for 
  there to be a corank-nullity expansion or even for the function to be 
  multiplicative or a graph invariant, that is, equal on isomorphic 
  graphs.  
   
  The motivation for allowing edge-specific values for the 
  deletion/con\-trac\-tion recursion comes from a number of applications 
  where it is natural.  This includes graphs with signed edges coming 
  from knot theory, graphs with edge-specific failure probability in 
  network reliability, and graphs whose edges represent various 
  interaction energies within a molecular lattice in statistical 
  mechanics.  While there is compelling motivation for allowing various edge 
  parameters,  
  the technical details of a general theory are 
  challenging. The two major works in this area are Zaslavsky \cite{Zas92} 
  and Bollob\'as and Riordan \cite{BR99}.  However, these two works take 
  different approaches, which were subsequently reconciled with a mild 
  generalization in \cite{E-MT06}, and for this reason we adopt the
  formalism   of \cite{E-MT06}. Bollob\'{a}s and Riordan~\cite{BR99} also give a
  succinct historical overview of the development of these
  multivariable extensions.  
   
  For the purposes of the following, we consider a class of graphs 
  minor-closed if it is closed under the deletion of loops, the 
  contraction of bridges, and the contraction and deletion of ordinary 
  edges; however we do 
  not require closure under the deletion of bridges.  Some formalism 
  is necessary to handle the parameters.  
   
  \begin{definition} 
  Let $U$ be a class, and let $R$ be a commutative ring. Then 
  an $R$-parametrization of $U$ consists of four parameter functions 
  $x,y,X,Y:U \to R$, denoted $e \to x_e ,y_e ,X_e ,Y_e $. 
  \end{definition} 
   
  \begin{definition} 
  Let $U$ be an $R$-parametrized class, and let $\Gamma$ be a 
  minor-closed class of graphs with $E(G) \subseteq U$ for 
  all $G \in \Gamma $.  Then a parametrized Tutte polynomial on 
  $\Gamma$ is a function $T:\Gamma  \to R$ which satisfies the 
  following: $T(G)  = X_eT(G/e)$ for any bridge $e$ of $G \in \Gamma $, 
  and $T(G) = Y_eT(G - e)$  for any loop $e$ of $G \in \Gamma $, and 
  $T(G) = y_eT(G - e) +  x_eT(G/e)$ for any ordinary edge $e$.   
   \end{definition} 
 
  The following theorem gives the central result.  Identity in 
  Item~\ref{itemA} comes  
  from requiring to be equal the two ways of carrying out 
  deletion/con\-trac\-tion reductions on a graph on two vertices with two 
  parallel edges $e_1$ and $e_2$ having parameters $ \left\{ {x_{e_i},
  y_{e_i } , X_{e_i } , Y_{e_i } } \right\}$.  Similarly, identities
  in Items~\ref{itemB} and~\ref{itemC} come from  
  considering the $\theta$-graph and $K_3$.  Here again $E_n$ is 
  the edgeless graph on $n$  vertices.  
   
  \begin{theorem}[The generalized Zaslavsky-Bollob\'as-Riordan
  theorem]\label{ZBR}  
  Let $R$ be a commutative ring, let $\Gamma $ be a 
  minor-closed class of graphs whose edge-sets are contained in an 
  $R$-parametrized class $U$, and let $a_1$, $a_2$,  \ldots  $\in R$.  Then 
  there is a parametrized Tutte polynomial   $T$ on $\Gamma $ with $T(E_n)  
  = a_n$  for all $n$ with $E_n  \in \Gamma$ if and only if the 
  following identities are satisfied:  
  \begin{enumerate} 
  \item Whenever $e_1$ and $e_2$ appear together in a circuit of a 
  $k$-component graph $G \in \Gamma $, then  
  $ a_k (x_{e_1} Y_{e_2} +  y_{e_1} X_{e_2}) = a_k  
  (x_{e_2} Y_{e_1} + y_{e_2} X_{e_1} )$.   
  \label{itemA} 
 
  \item Whenever $e_1$, $e_2$ and $e_3$ appear together in a circuit of a 
  $k$-component graph $G \in \Gamma $, then  
  $ a_k X_{e_3}  
  (x_{e_1} Y_{e_2} + y_{e_1}x_{e_2}) = a_k X_{e_3}  
  (Y_{e_1} x_{e_2} + x_{e_1} y_{e_2} )$.  \label{itemB} 
 
  \item Whenever $e_1$, $e_2$ and $e_3$ are parallel to each other in a 
  $k$-component graph $G \in \Gamma $, then  
   $ a_k Y_{e_3} 
  (x_{e_1} Y_{e_2} +  y_{e_1} x_{e_2}) = a_k  Y_{e_3}  
  (Y_{e_1} x_{e_2} +  x_{e_1} y_{e_2} )$.  \label{itemC} 
    \end{enumerate} 
  \end{theorem}

  A most general parametrized Tutte polynomial, what possibly could be 
  called {\em the} parametrized Tutte polynomial, might begin with the 
  polynomial ring on independent variables $\left\{ {x_e ,y_e ,X_e ,Y_e 
      {\text{: }}e \in U} \right\} \cup \left\{ {a_i :i \geqslant 1} 
  \right\}$.  However, the resulting function is not technically a 
  polynomial, in that it must take its values not in the polynomial 
  ring, but has as $R$ the polynomial ring modulo the ideal generated by 
  the identities in Theorem~\ref{ZBR}.   
 
The question also arises as to 
  whether ``the most general''  parametrized Tutte polynomial should be 
  multiplicative on disjoint unions and one point joint of graphs, as 
  this introduces additional relations among the $a_i$'s.  This is 
  because a parametrized Tutte polynomial is not necessarily 
  multiplicative.  A sufficient condition is the following. 
  
\begin{proposition} 
Suppose $T$ 
  is a parametrized Tutte polynomial on a minor-closed class of graphs 
  that contains at least one graph with $k$ components for every $k$ and 
  that is closed under  one-point unions and the removal of isolated 
  vertices. Then $T$ is multiplicative with respect to both disjoint 
  unions and one-point joins if and only if the $\alpha_1 = T(E_1)$ is 
  idempotent, and $\alpha_k = \alpha_1$ for all $k \geqslant 1$.  
\end{proposition} 
   
  Bollob\'as and Riordan \cite{BR99} emphasize graph invariants, and hence 
  require that the parametrization be a coloring of the graph.  That 
  is, graphs are edge-colored (not necessarily 
  properly), with edges of the same color having the same parameter 
  sets.  This enables consideration of parametrized Tutte polynomials
   that are invariants of colored graphs, but requires the following 
  additional constraints.  For every $e_1 \in U$, there are $e_2$, $e_3 
  \in U$  with $e_1  \ne e_2  \ne e_3  \ne e_1 $ such that $x_{e_1} = 
  x_{e_2} = x_{e_3}$ , $y_{e_1} = y_{e_2} =y_{e_3}$ , $X_{e_1} = 
  X_{e_2} = X_{e_3}$, and   $Y_{e_1} = Y_{e_2} = Y_{e_3}$.

 Proofs of the above results and further details may be found it
 \cite{Zas92, BR99, E-MT06}.  We note that any relation between this
 Tutte polynomial generalization with its edge parameters, and the W-
 and U-polynomials of Subsection \ref{subsec:symmetric} with their
 vertex weights, has not yet been studied.   
 
Interestingly, although the parametrized Tutte polynomial has an
activities expansion analogous to that of the classical Tutte
polynomial, it does not necessarily have an analog of the rank-nullity
formulation.  However, under modest assumptions involving non-zero
parameters and some inverses, the parametrized Tutte polynomial may be
expressed in a rank-nullity form.  This is fortunate, because
significant results for the zeros of the chromatic and Tutte
polynomial have arisen from such a multivariable realization.
Examples may be found in Sokal~\cite{Sok01a}, Royle and
Sokal~\cite{RS04}, and Choe, Oxley, Sokal, and Wagner~\cite{COSW04}.

\subsection{The Generalized Transition Polynomial}
\label{gen_transition_poly}
 
A number of state model polynomials, for example the circuit partition
polynomials, Penrose polynomial, the Kauffman bracket for knots and
links, and the transition polynomials of Jaeger \cite{Jae90}, that are
not specializations of the Tutte polynomial, are specializations of
the multivariable generalized transition polynomial of \cite{E-MS02}
which we describe here.  This multivariable extension is a Hopf
algebra map, which leads to structural identities that then inform its
various specializations.  The medial graph construction that relates
the circuit partition polynomial and the classical Tutte polynomial
\index{Tutte polynomial}\index{polynomial!Tutte}  extends to similarly
relate the generalized transition polynomial and the parametrized
Tutte polynomial when it has a rank-nullity formulation.  
 
The graphs  here are Eulerian, \index{Eulerian graph}
\index{graph!Eulerian} although not necessarily connected,
with loops and multiple edges allowed.  A vertex state, 
\index{vertex state} or transition, is a choice of local
reconfiguration of a 
graph at a vertex by pairing the half edges incident with that vertex.
A graph state, \index{graph state} or transition system, $S\left( G
\right)$, is the result of choosing a vertex state at each vertex of
degree greater than 2, and hence is a union of disjoint cycles. 
 We will write $St(G)$ for the set of graph states of $G$, and throughout 
we assume weights have values in $R$, a commutative
ring with unit.  
 
A \emph{skein relation} \index{skein relation} for graphs is a formal sum of
weighted vertex states, together with an evaluation of the terminal
forms (the graph states).  See~\cite{E-M98, E-MS02} for a detailed
discussion of these 
concepts, which are appropriated from knot theory, in their most
general form, and Yetter~\cite{Yet90} for a general theory of
invariants given by linear recursion relations. A \emph{skein type}, (or
state model, or transition) \emph{polynomial} is one which is computed by
repeated applications of skein relations.  See Jaeger~\cite{Jae90}
for a comprehensive treatment of these in the case of 4-regular
graphs.  
 
For brevity, we elide technical details such as free loops and
isomorphism classes of graphs with weight systems which may be found
in \cite{E-MS02}.  
 
\begin{definition} 
Pair, vertex, and state weights: 
\begin{enumerate} 
\item{A \emph{pair weight}\index{pair weight} is an association of a
    value $p\left( {e_v ,e'_v } \right)$  in a unitary ring $R$
    to a pair of half edges incident with a vertex $v$ in $G$.   A
    \emph{weight system},\index{weight system} $W( G )$, of
    an Eulerian graph $G$ is an assignment of a pair weight to every
    possible pair of adjacent half edges of $G$.}

\item{The \emph{vertex state weight}  \index{vertex state weight} of a
    vertex state is $\prod {p( {e_v ,e'_v } )}$
    where the product is over the pairs of half edges comprising
    the vertex state.}

\item{The \emph{state weight} \index{state weight} of a graph state
    $S$ of a graph $G$ with weight system $W$ is $\omega( S ) =
    \prod{\omega( {v,S} )}$,
    where $\omega( {v,S} )$ is the vertex state weight of
    the vertex state at $v$ in the graph state $S$, and where the
    product is over all vertices of  $G$.}  
\end{enumerate} 
\end{definition}

When $A$ is an Eulerian subgraph of an Eulerian graph $G$ with weight
system $W(G)$, then $A$ inherits its weight system $W(A)$ from $G$ in
the obvious way, with each pair of adjacent edges in $A$ having the
same pair weight as it has in $G$.   When $A$ is a graph resulting
from  locally replacing the vertex $v$ by one of its
$\prod_{i = 0}^{n - 1} {( {2n - (2i + 1)} )} $
vertex states, then all the pair weights are the same as they are in
$G$, except that all the pairs of half edges adjacent to the newly
formed vertices of degree 2 in $A$ have pair weight equal to 1, the
identity in $R$.  
 
The generalized transition polynomial $N\left( {G; W, x} \right)$ has
several formulation, and we give two of them, a linear recursion
formula and a state model formula, here.
\index{generalized transition polynomial}
\index{transition polynomial}
\index{polynomial!generalized transition}
\index{polynomial!transition}  
 
\begin{definition} \label{Recursive transition} 
The generalized transition polynomial, $N( {G;W,x} )$, is defined
recursively by repeatedly applying 
the skein relation 
 \[ N( {G;W,x} ) =\sum {\beta _i N( {G_i ;W(G_i),x})}\]  
at any vertex $v$ of degree greater than 2.  Here the $G_i $'s are the
graphs that result from locally replacing a vertex $v$ of degree $2n$
in G  by one of its  vertex states.  The $\beta_i $'s are the vertex
state weights.  Repeated application of this relation reduces $G$ to a
weighted (formal) sum of disjoint unions of cycles,  (the graph
states).  These 
terminal forms are evaluated by identifying each cycle with the
variable $x$, weighted by the product of the pair weights over all
pairs of half edges in the cycle.   
\end{definition}

\begin{definition} \label{State Model transition}  The state model
  definition of the generalized transition polynomial is:  
\[ 
N(G;W,x) = \sum\limits_{St\left( G \right)}^{} {\left( {\left(
        {\prod\limits_{} {\omega \left( {v,S} \right)} }
      \right)x^{k(S)} } \right)}  = \sum\limits_{St\left( G
  \right)}^{} {\omega \left( S \right)x^{k(S)} }. \]  
 
\end {definition} 
 
Note that vertex states commute, that is, if $G_{uv}$ results from
choosing a vertex state at $u$, and then at $v$, we have
$G_{uv}=G_{vu}$.  Thus, Definition \ref{Recursive transition} gives a
well-defined function, and Definitions \ref{Recursive transition} and
\ref{State Model transition} are equivalent.

Several of the polynomials we have already seen are specializations of
this generalized transition polynomial.  For example, if all the pair
weights are 1, then the circuit partition polynomial
 \index{circuit partition polynomial}
 \index{polynomial!circuit partition} for an
unoriented Eulerian graph results.  If $\vec G$ is an Eulerian
digraph, and $G$ is the underlying undirected graph with pair weights
of 1 for pairs half edges corresponding to one inward and one outward
oriented half edge of $\vec G$ and 0 otherwise, then the oriented
version of the circuit partition polynomial results.    
 
In the special case that $G$ is 4-regular, the polynomial $N\left(
 {G;W,x} \right)$  
 is essentially the same as the transition polynomial
 \index{transition polynomial}\index{polynomial!transition}
 $Q(G,A,\tau )$ of Jaeger~\cite{Jae90}, where $G$ is a 4-regular
 graph and $A$ is a system of vertex state weights (rather than pair
 weights).  If the vertex state weight in $\left( {G,A} \right)$  
 is $w$, then define $W(G)$ by letting the pair weights for each of
 the two pairs of edges determined by the state be $\sqrt w $.  The
 two polynomials then just differ by a factor of $x$, so $N\left(
 {G;W,x} \right)=xQ(G,A,x)$, and here we retain vertices of degree 2 in
 the recursion while they are elided in \cite{Jae90}.  Thus $N\left(
 {G;W,x} \right)$ gives a generalization of Jaeger's transition
 polynomials to all Eulerian graphs.  
 
Because $Q(G,A,\tau )$ assimilates the original Martin polynomial for 4-regular graphs and digraphs, \index{Martin polynomial}
\index{polynomial!Martin} the Penrose polynomial \index{Penrose polynomial}
\index{polynomial!Penrose} and the Kauffman
bracket of knot theory 
\index{Kauffman bracket}
\index{polynomial!Kauffman bracket} (see \cite{Jae90}), and
$N\left( {G;W,x} \right)$ assimilates $Q(G,A,\tau )$, we have that the
Penrose polynomial and Kauffman bracket are also specializations of
$N\left( {G;W,x} \right)$.  Specifically, if $G$ is a planar graph
with face 2-colored medial graph $G_m$, and we give a weight system to
$G_m$ by assigning a value of 1 to pairs of edges that either cross at
a vertex or bound the same black face and 0 otherwise, then $N\left(
  {G_m ;W,x} \right) = P(G;x)$.  Similarly, if $L$ is a link, and
$G_L$ is the signed, face 2-colored universe of $L$, then a weight
system can be assigned to $G_L$ so that $N\left( {G_L ;W,a^2  + a^{ -
      2} } \right) = (a^2  + a^{ - 2} )K[L]$ where $K[L]$ is the
Kauffman bracket of the link.   
 
Because of these specializations, the algebraic properties of the
generalized transition polynomial are available to inform these other
polynomials as well. In particular,  $N\left( {G;W,x} \right)$ is a
Hopf algebra map from  the freely
generated (commutative) hereditary Hopf algebra of Eulerian graphs
with weight systems to the binomial bialgebra $R[x]$ (details may be found in \cite{E-MS02}).  This
leads to two structural identities, the first from the
comultiplication in the Hopf algebra, the second from the antipode.      
 
\begin{theorem} 
Let $G$ be an Eulerian graph.  Then 
\[ 
N\left( {G;W,x + y} \right) = \sum\limits_{}^{} {N\left( {A_1 ;W\left(
        {A_1 } \right),x} \right)N\left( {A_2 ;W\left( {A_2 }
      \right),y} \right)}   
\] 
where the sum is over all ordered partitions of $G$ into two
edge-disjoint Eulerian subgraphs $A_1$ and $A_2$, and  
\[ 
N\left( {G;W, - x} \right) = N\left( {\zeta \left( {G;W} \right),x} \right), 
\] 
where $\zeta$ is the antipode $\zeta \left( {G,W} \right) = \sum {( -
  1)^{|P|}( A_1  \ldots A_{|P|})} $, with the sum over all ordered
partition $P$ of $G$ into $|P|$ edge-disjoint Eulerian subgraphs each
with inherited weight system. Here $N(G;W,x)$ is extended linearly over such formal sums. 
\end{theorem}

This type of Hopf algebraic structure has already been used to
considerably extend the known combinatorial interpretations for
evaluations of the Martin, Penrose, and Tutte polynomials   implicitly
by Bollob\'as~\cite{Bol02}, and explicitly by Ellis-Monaghan and
Sarmiento~\cite{E-M98,  E-MS01,  Sar01,  E-M04a,  E-M04b}. The first
identity has been used to find combinatorial interpretations for the
Martin polynomials for all integers, where this was previously only
known for –2, -1, 0, 1 in the oriented case, and –2, 0, 2 in the
unoriented case. This then led to combinatorial interpretations for
the Tutte polynomial (and its derivatives) of a planar graph for all
integers along the line $x=y$, where previously –1, 3 were the only
known non-trivial values. These results for the Tutte polynomial were
mentioned in the preceding chapter and for the circuit partition
polynomial in Subsection \ref{sec:cpp}. The second identity has been
used to determine combinatorial interpretations for the Penrose
polynomial for all negative integers, where this was previously only
known for positive integers.

\section{Two Applications}
\label{sec:applications}

Graph polynomials have a wide range of applications throughout many
fields.  We have already seen some examples of this with various
applications of the classical Tutte polynomial in the previous
chapter.  Here we present two representative important applications
(out of many possible) and show how they may be modeled by graph
polynomials. 

\subsection{DNA Sequencing}
\label{subsec:dna}

We begin with string reconstruction, a problem that may be modeled by
the interlace and circuit partition polynomials 
\index{interlace polynomial}
\index{polynomial!interlace}
\index{circuit partition polynomial}
\index{polynomial!circuit partition} (and hence
indirectly in special cases by the Tutte polynomial).
String reconstruction is the 
process of reassembling a long string of symbols from a set of its
subsequences together with some (possibly incomplete, redundant, or
corrupt) sequencing information.  While we focus on DNA sequencing,
which was original the motivation for the development of the interlace
polynomial, the methods here apply to any string reconstruction
problem. For example, fragmenting and reassembling messages is a
common network protocol, and reconstruction techniques might be
applied when the network protocol has been disrupted, yet the original
message must be reassembled from the fragments. 

DNA sequences are typically too long to read at once with current laboratory techniques, so researchers
probe for shorter fragments (reads) of the strand.  They then are
faced with the difficulty of recovering the original long sequence
from the resulting set of subsequences.  DNA sequencing by
hybridization is a method of reconstructing the nucleotide sequence
from a set of short substrings (see Waterman~\cite{Wat95} for an
overview).  
The problem of determining the number of possible reconstructions may
be modeled using Eulerian digraphs, with a correct sequencing of the
original strand corresponding to exactly one of the possible Eulerian
circuits in the graph.  The probability of correctly sequencing the
original strand is thus the reciprocal of the total number of Euler
circuits in the graph.  
 
 The most basic (two-way repeats only) combinatorial model  for DNA
 sequencing by hybridization uses an Eulerian digraph with two
 incoming and two outgoing edges at each vertex (see
 Pevzner~\cite{Pev89} and Arratia, Martin, Reinert,
 Waterman~\cite{AMRW96}). The raw data consists of all subsequences of
 the DNA strand of a fixed length L, called the L-spectrum of the
 sequence.  As L increases, the statistical probability that the
 beginning and end of the DNA strand are the same approaches zero, as
 does the likelihood of three or more repeats of the same pattern of
 length L or more in the strand (see Dyer, Frieze, and
 Suen~\cite{DFS94}).  Thus, this model assumes that the only
 consideration in reconstructing the original sequence is the
 appearance of interlaced two-way repeats, that is, alternating
 patterns of length L or greater, for example, $\dots \text{ACTG}
 \dots \text{CTCT} \dots \text{ACTG} \dots \text{CTCT} \dots$ . 

From the multiset (duplicates are allowed) of subsequences of length
L, create a single vertex of the de Bruijn graph for each subsequence
of length L-1 that appears in one of the subsequences.  For example,
if L = 4 and ACTG appears as a subsequence, create two vertices, one
labeled ACT and one labeled CTG.  Edges are directed from head to tail
of a subsequence, e.g. there would be a directed edge from ACT to CTG
labeled ACTG.  If there is another subsequence ACTT, we do not create
another vertex ACT, but rather draw an edge labeled ACTT from the
vertex ACT to a new vertex labeled CTT.  If, in the multiset of
subsequences, ACTG appears twice, then we draw two edges from ACT to
CTG. 

The beginning and end of the strand are identified to be represented
by the same vertex, and, since by assumption no subsequence appears
more than twice, the result is an Eulerian digraph of maximum degree
4.  Tracing the original DNA sequence in this graph corresponds to an
Eulerian circuit that starts at the vertex representing the beginning
and end of the strand.  All other possible sequences that could be
(mis)reconstructed from the multiset of subsequences correspond to
other Eulerian circuits in this graph.  Thus (up to minor reductions
for long repeats and forced subsequences), finding the number of DNA
sequences possible from a given multiset of subsequences corresponds
to enumerating the Eulerian circuits in this directed
graph. 
 
The generalized transition polynomial models this problem directly:
when the pair weights are identically 1, it reduces to the circuit
partition polynomial.  This is a generating function for families of
circuits in a graph, so the coefficient of $x$ is the number of
Eulerian circuits. The interlace polynomial informs the problem as
follows. Consider an Eulerian circuit through the de Bruijn graph,
which gives a sequence of the vertices visited in order.  
Now construct the interlace graph by placing a vertex for each symbol
and an edge between symbols that are interlaced (occur in alternation)
in the sequence. The interlace polynomial of the interlace graph is
then a translation of the circuit partition polynomial of the original
de Bruijn graph, as in Theorem~\ref{A}, where again the coefficient of $x$ is the number of
Eulerian circuits (see Arratia, Bollob\'as, Coppersmith, and
Sorkin~\cite{ABCS00}, and Arratia, Bollob\'as, and
Coppersmith~\cite{ABS00, ABS04a}).

One of the original motivating goals of Arratia, Bollob\'as,
Coppersmith, and Sorkin ~\cite{ABCS00} was classifying Eulerian
digraphs with a given number of Eulerian circuits. The BEST theorem,
a formula for the number of the circuits of an Eulerian graph in terms
of its Kirchhoff matrix (see Fleischner~\cite{Fle91} for good
exposition) gives only a tautological classification: the Eulerian
digraphs with $m$ Eulerian circuits are those where BEST theorem
formula gives $m$ circuits.  Critically, all of the above graph
polynomials encode much more information than is available from the
BEST theorem, and all of them are embedded in broader algebraic
structures that provide tools for extracting information from them.
Thus, they better serve the goal of seeking structural
characterizations of graph classes with specified Eulerian circuit
properties.

\subsection{The Potts Model of Statistical Mechanics}
\label{subsec:potts}

Here we have an important physics model that remarkably was found to
be exactly equivalent to the Tutte polynomial.
 \index{Tutte polynomial}\index{polynomial!Tutte}

Complex systems are networks in which very simple interactions at the
microscale level determine the macroscale behavior of the system. The
Potts model \index{Potts model partition function} of statistical
mechanics models complex systems whose behaviors depend on nearest
neighbor energy interactions.  This model plays an important role in
the theory of phase transitions and critical phenomena, and has
applications as widely varied as magnetism, adsorption of gases on
substrates, foam behaviors, and social demographics, with important
biological examples including disease transmission, cell migration,
tissue engulfment, diffusion across a membrane, and cell sorting. 

Central to the Potts model is the Hamiltonian, \index{Hamiltonian}
\[
 h(\omega ) =  - J\sum\limits_{\{ i,j\}  \in E(G)}
              {\delta (\sigma _i,\sigma _j )} \;,
\]
 a measure of the energy of the system.  Here a
spin, $\sigma _i $, at a vertex $i$, is a choice of condition (for
example, healthy, infected or necrotic for a cell represented by the
vertex). $J$ is a measure of the interaction energy between
neighboring vertices, $\omega $ is a state of a graph $G$ (that is, a
fixed choice of spin at each vertex), and $\delta $ is the Kronecker
delta function. 

The Potts model partition function is the normalization factor for the
Boltzmann probability distribution.  Systems such as the Potts model,
following Boltzmann distribution laws, will have the number of states
with a given energy (Hamiltonian value) exponentially distributed.
Thus, the probability of the system being in a particular state
$\omega$ at temperature $t$ is: 
\[ 
\Pr \left( {\omega ,\beta } \right) = \frac{{\exp \left( { - \beta
        h\left( \omega  \right)} \right)}}{{\sum {\exp \left( { -
          \beta h\left( \varpi  \right)} \right)} }}.
\]
Here, the sum is over all possible states $\varpi$ of G, and $\beta  =
\frac{1}{{\kappa t}}$ 
, where $\kappa  = 1.38 \times 10^{ - 23} $ joules/Kelvin is the
Boltzmann constant. The parameter $t$ is an important variable in the
model, although it may not represent physical temperature, but some
other measure of volatility relevant to the particular application
(for example ease of disease transmission/reinfection).  The
denominator of this expression, $P\left( {G;q,\beta } \right) = \sum
{\exp \left( { - \beta h\left( \varpi  \right)} \right)} $, called the
Potts model partition function, is the most critical, and difficult,
part of the model. 

Remarkably, the Potts model partition function is equivalent to the
Tutte polynomial: \index{Tutte polynomial}\index{polynomial!Tutte}  
\[
P\left( {G;q,\beta } \right) = q^{k\left( G \right)} v^{\left|
    {v\left( G \right)} \right| - k\left( G \right)} T\left(
  {G;\frac{{q + v}}{v},v + 1} \right), 
\]
where $q$ is the number of possible spins, and $v = \exp (J\beta ) -
1$.  See Fortuin and Kasteleyn \cite{FK72} for the nascent stages of
this theory, later exposition in Tutte ~\cite{Tut84}, Biggs
~\cite{Big96}, Bollob\'as ~\cite{Bol98}, Welsh ~\cite{Wel93}, and
surveys by Welsh and Merino ~\cite{WM00}, and Beaudin, Ellis-Monaghan,
Pangborn and Shrock ~\cite{BE-MPS}.  

One common extension of the Potts model involves allowing interaction
energies to depend on individual edges.  With this, the Hamiltonian
becomes $h(\omega ) = \sum_{e \in E(G)} {J_e \delta (\sigma _i
  ,\sigma _j )} $, where $J_e $  is the interaction energy on the edge
$e$.  The partition function is then  
 \[
P\left( G \right) = \sum\limits_{A \subseteq E\left( G \right)}
{q^{k\left( A \right)} \prod\limits_{e \in A} {v_e } },  
\]
where $v_e  = \exp \left( {\beta J_e } \right) - 1$
 Again see Fortuin and Kasteleyn ~\cite{FK72}, and more recently Sokal
 ~\cite{Sok00, Sok01b}.  As we have seen in Subsection
 \ref{subsec:param tutte}, the Tutte polynomial has also been extended
 to parametrized Tutte functions that incorporate edge weights.   The
 generalized partition function given above satisfies the relations of Theorem~\ref{ZBR}, however, and thus is a special case of a parametrized
 Tutte function. 

This relationship between the Potts model partition function and the
Tutte polynomial \index{Tutte polynomial}\index{polynomial!Tutte}  has
led to a remarkable synergy between the fields, particularly for
example in the areas of computational complexity and the zeros of the
Tutte and chromatic polynomials.  For  overviews, see Welsh and
Merino ~\cite{WM00}, and Beaudin, Ellis-Monaghan, Pangborn and
Shrock~\cite{BE-MPS}.

\section{Conclusion}
\label{conclusion}

There are a great many other graph polynomials equally interesting to
those surveyed here, including for example the F-polynomials
\index{F-polynomial}\index{polynomial!F} of Farrell, the Hosaya or
Wiener polynomial, the clique/ independence and adjoint polynomials,
etc.  In particular, Farrell~\cite{Far79b} has a circuit cover
polynomial (different from the circuit partition polynomial of
Subsection~\ref{sec:cpp}) with noteworthy interrelations with the
characteristic polynomial. Also, Chung and Graham developed a
`Tutte-like' polynomial for directed graphs in \cite{CG95}.  The
resultant cover polynomial \index{cover polynomial}
\index{polynomial!cover}is extended to a symmetric
function generalization, like those in section \ref{subsec:symmetric},
by Chow \cite{Cho96}. Similarly, Courcelle~\cite{Cou08} and
Traldi~\cite{Tra} have also very recently developed multivariable
extensions of the interlace polynomial.  Some surveys of graph
polynomials with 
complementary coverage to this one include: Pathasarthy ~\cite{Par89},
Jaeger ~\cite{Jae90}, Farrell ~\cite{Far93}, Fiol ~\cite{Fio97},
Godsil ~\cite{God84}, Aigner ~\cite{Aig01}, Noy ~\cite{Noy03} and Levit
and Mandrescu ~\cite{LM05}.

%
%
%

%
%



\printindex
\end{document}